\documentclass[11pt,twoside,a4paper]{article}
\usepackage{color,amscd,amsmath,amssymb,amsthm,latexsym,stmaryrd,authblk,mathabx,shuffle,hyperref,tikz,tkz-tab} 
\usepackage[latin1]{inputenc}
\usepackage[T1]{fontenc}   
\usepackage[english]{babel}
\input{xy}
\xyoption{all}

\title{Matrix symmetric and quasi-symmetric functions and noncommutative representation theory}
\date{}
\author{Lo\"\i c Foissy, Claudia Malvenuto, Fr\'ed\'eric Patras}

\theoremstyle{plain}
\newtheorem{theo}{Theorem}[section]
\newtheorem{lemma}[theo]{Lemma}
\newtheorem{cor}[theo]{Corollary}
\newtheorem{prop}[theo]{Proposition}
\newtheorem{defi}[theo]{Definition}

\theoremstyle{remark}
\newtheorem{remark}{Remark}[section]
\newtheorem{notation}{Notations}[section]
\newtheorem{example}{Example}[section]

\renewcommand{\leq}{\leqslant}
\renewcommand{\geq}{\geqslant}
\newcommand{\K}{\mathbb{K}}
\newcommand{\N}{\mathbb{N}}
\newcommand{\Q}{\mathbb{Q}}

\newcommand{\bfH}{\mathbf{H}}

\newcommand{\id}{\mathrm{Id}}

\renewcommand{\ker}{\mathrm{Ker}}
\newcommand{\prim}{\mathrm{Prim}}

\newcommand{\pack}{\mathrm{Pack}}
\newcommand{\mat}{\mathrm{Mat}}
\newcommand{\bfA}{\mathbf{A}}
\newcommand{\bfM}{\mathbf{M}}
\newcommand{\calM}{\mathcal{M}}
\newcommand{\qsh}{\mathrm{qsh}}
\newcommand{\adm}{\mathrm{Adm}}
\newcommand{\QSym}{\mathbf{QSym}}
\newcommand{\NSym}{\mathbf{NSym}}
\newcommand{\row}{\mathrm{row}}
\newcommand{\col}{\mathrm{col}}

\newcommand{\comp}{\mathrm{Comp}}
\newcommand{\sh}{\mathrm{sh}}
\newcommand{\inc}{\mathrm{inc}}
\newcommand{\Exp}{\mathrm{Exp}}

\newcommand{\bdeltares}{\blacktriangle_{\mathrm{res}}}
\newcommand{\desc}{\mathit{Desc}}
\newcommand{\gr}{\mathrm{gr}}
\begin{document}

\maketitle

\begin{abstract}
A fundamental result by L. Solomon in algebraic combinatorics and representation theory states that Mackey formulas for products of characters of a symmetric group, or equivalently the computation of tensor products of representations thereof, can be lifted to the corresponding Solomon's descent algebra, a subalgebra of the group algebra with a very rich structure. Motivated by the structure of the product formula in these algebras and by other results and ideas in the field, we introduce and investigate in the present article a two dimensional analog of descent algebras based on packed integer matrices that inherits most of their fundamental properties.
One of the various bialgebra structures we introduce on packed integer matrices identifies with a bialgebra recently introduced by 
J. Diehl and L. Schmitz to define a two dimensional generalisation of Chen's iterated integrals signatures.

\end{abstract}

\tableofcontents

\section{Introduction}
People familiar with the combinatorial approach to symmetric group representations and with the theory of quasi-symmetric functions know that many objects involved in the theory are two dimensional. This is the case for example of Young diagrams (parametrising irreducible representations and parabolic subgroups), of standard and semi-standard tableaux (used to effectively construct these representations and perform explicit computations). At a deeper level of analysis, Zelevinski's theory of pictures also features the representation theoretic role of pair of orders on finite sets. These considerations have originated some of our previous works on plane posets, double posets and pictures such as \cite{Malvenuto2011,Foissy20,Foissy33}. Although addressing a seemingly different topic, the present article should be understood as a contribution 
in the same overall direction.

A fundamental result by L. Solomon in algebraic combinatorics and representation theory \cite{Solomon1976} states that Mackey formulas for products of characters of a symmetric group $S_n$ can be lifted to the corresponding Solomon's descent algebra, a noncommutative subalgebra of the group algebra with a very rich structure. The direct sum of all Solomon's descent algebras carries for example a graded cocommutative Hopf algebra structure with compatibility relations with products in the symmetric group algebras: its dual is a double bialgebra in the sense of \cite{Foissy40} (see also Section \ref{sectiondouble} in the present article). On descent algebras, see e.g. \cite{Malvenuto1995} and \cite[Chapter 5]{Cartier2021}. 
This problem of computing products of characters is one of the most classical in algebraic combinatorics. It is equivalent to the one of computing tensor products of symmetric group representations or products of symmetric functions. 
For example, 
the Littlewood--Richardson rule computes the product of two Schur functions in the basis of Schur functions by enumerating skew tableaux: an example of 2-dimensional calculus to which the theory of pictures provides an alternative approach \cite{Malvenuto2011}.

The product formula in Solomon's algebras is thus a very interesting object in itself. Compositions of integers naturally index bases of these algebras. In such a basis, products are naturally parametrised by packed integer matrices (as briefly recalled in the next section of the present article) and involve the mapping from matrices to compositions obtained from the reading order.
In the present article we introduce a two dimensional analog of descent algebras based on packed integer matrices, many properties of which seem to reflect combinatorial properties of this product. The graded vector space spanned by packed integer matrices, its graded components and the dual spaces have a very rich structure that naturally generalise those on the Solomon's descent algebras and their direct sum.
Among various other results of the same type, we show for example that the map from packed matrices to compositions induced by the reading order is a bialgebra map (Prop. \ref{propTheta}).

Lastly, recall that the polynomial realisation of the dual of the direct sum of the various Solomon's descent algebras is isomorphic to the algebra of quasi-symmetric functions \cite{Malvenuto1995}. The isomorphism is actually an isomorphism of double bialgebras \cite{Foissy40}.
In the same vein we study polynomial realisations of packed integer matrices and define and investigate various associated double bialgebra morphisms.

The idea of constructing combinatorial Hopf algebras using packed matrices is not new. It seems to have first appeared in an article by G. Duchamp, F. Hivert and J.-Y. Thibon \cite{Duchamp2002}, together with a notion of matrix quasi-symmetric functions. However the bialgebras they defined have a different structure than the present ones: in particular the treatment of rows and columns is not symmetric in their constructions that do not seem to be related to ours. 

Our results instead directly connect to recent work by J. Diehl and L. Schmitz \cite{Diehl}. They defined a graded bialgebra structure on packed matrices that identifies with one of the bialgebra structures on packed matrices we consider (the graded commutative one of Section \ref{ssect:gdhp}). We also use the same polynomial realisation of packed matrices (Section \ref{ssec:const}). Their work aims at defining a two dimensional generalisation of Chen's iterated integrals signatures for image processing. Our results complement theirs in many respects and we hope they can be useful from this point of view.

Still another graded commutative bialgebra structure on packed matrices has been introduced by Giusti et al., the two parameters shuffle bialgebra \cite{Giusti,Schmitz}. L. Schmitz and N. Tapia have proven that it is isomorphic to the bialgebra defined by Diehl and Schmitz \cite{Schmitz}. We give a brief account of these results in Section \ref{ssec:sprod}, where we also
show that a cocommutative bialgebra of packed matrices
constructed at the end of Section \ref{sec:bim} is a bigraded dual of the two parameters shuffle bialgebra. \\

The article is organised as follows. 

Section \ref{descents} summarizes general properties of descents in symmetric groups and their relation to $\NSym$, the bialgebra of noncommutative symmetric functions, as a motivation for later developments on their two dimensional generalisation. The Section also contains reminders on double bialgebras, a key notion in the theory of descents in symmetric groups and quasi-symmetric functions, and also in our article. We also take into account  classical variants of the main constructions that will be generalized later on in the article to the matrix case.

Section \ref{sec:bim} introduces two dimensional matrix and packed matrix analogs of the graded cocommutative bialgebra of descents and noncommutative symmetric functions. Bialgebra maps to and from $\NSym$ are introduced. We also discuss several variants of the main constructions, extending in particular to the matrix case the constructions in section \ref{descents}.

Section \ref{sec:dualty} studies the dual structures and a natural infinitesimal bialgebra structure on packed matrices. Links with the bialgebra of matrix compositions of Diehl-Schmitz and with the Giusti-Lee-Nanda-Oberhauser two parameters shuffle Hopf algebra are explicited.  Bialgebra maps from and to the bialgebra of quasi-symmetric functions are also obtained.

Section \ref{sec:double} focuses on some double bialgebra structures on packed matrices. The section starts with polynomial realizations, recalling in particular the Diehl-Schmitz construction of a bialgebra of two parameters symmetric functions. Another, new, bialgebra structure and a double bialgebra structure on packed matrices are obtained. The section and the article conclude with the proof that morphisms previously introduced are actually preserving the full algebraic structure involved: they are morphisms of double bialgebras.\\

We assume that the reader is familiar with the language and fundamental notions of the theory of coalgebras, bialgebras and Hopf algebras. A general reference for the notions and results used throughout the article is \cite{Cartier2021}.\\

\begin{notation}\ \
\begin{enumerate}
\item Let $k\in \N$.  We denote by $[k]$ the set $\{1,\ldots,k\}$. By convention, $[0]=\emptyset$.  
\item The $n$--th symmetric group, that is the group of permutations of $[n]$,  is denoted $S_n$. We write $1_n$ its unit (the identity map). Its group algebra over $\Q$ (i.e. the $\Q$-vector space spanned by its elements equipped with the product induced by the group law) is denoted by $\Q[S_n]$.
\item The convolution product of linear endomorphisms of a bialgebra $B$ with product $m$, coproduct $\Delta$, unit $\eta$ and counit $\varepsilon$ is denoted $\ast$: 
$$f\ast g:=m\circ (f\otimes g)\circ \Delta.$$
It is associative, with unit denoted by $\nu:=\eta\circ\varepsilon$. When the bialgebra is graded connected (that is, its degree 0 component identifies with the ground field), $\nu$ is the projection on the ground field orthogonally to the higher degree components. In that case, the bialgebra is automatically a Hopf algebra: its identity map has a convolution inverse.
\item We use sometimes the Sweedler notation and use $x^{(1)}\otimes x^{(2)}$ as a shortcut for $\Delta(x)$, for $x$ an arbitrary element in a coalgebra with coproduct $\Delta$.
\item A sequence $\mu=(\mu_1,\dots,\mu_k)$ of positive integers is called a composition of the integer $|\mu|:=\mu_1+\dots +\mu_k$. A sequence $\mu=(\mu_1,\dots,\mu_k)$ of nonnegative integers is called a generalised composition of $|\mu|$ (the  weight  of $\mu$).
\item
For any $n\geq 0$, we denote by $H_n(X)$ the $n$-th Hilbert polynomial:
\[H_n(X)=\frac{X(X-1)\ldots (X-n+1)}{n!}.\]
In particular, $H_0(X)=1$. 
\end{enumerate}
\end{notation}

\section{Reminders on descents and double bialgebras}\label{descents}
\subsection{Descents and noncommutative symmetric functions}
We provide here a brief account of the algebraic theory of descents in symmetric groups from the point of view of the theory of words and Hopf algebras, following \cite{Malvenuto1995,Cartier2021}. As mentioned in the Introduction, another approach would be possible, namely the one of Solomon that uses descents to construct a lift to the group algebra $\Q[S_n]$ of calculations in the ring of characters of $S_n$ \cite{Solomon1976} (or, equivalently, in the Grothendieck ring of linear representations of $S_n$). However, it would take longer to describe, and we stick to the word approach.

Recall that a permutation $\sigma$ in $S_n$ has a descent in position $i<n$ if and only if $\sigma(i)>\sigma(i+1)$. The set of descents of $\sigma$ is written $\desc(\sigma)\subseteq [n-1]$. Let $U\subseteq [n-1]$: the set of all permutations in $S_n$ whose descent set is $U$ is written $D_{=U}$ and we call it the descent class associated to $U$. The set of all permutations whose descent set is included in $U$ is written $D_{\subseteq U}$.

Consider now the shuffle Hopf algebra $(T(X),\shuffle,\Delta)$, where $X$ is a countable infinite alphabet; $T(X)$ is the graded vector space spanned by words  $y_1\dots y_n$, $y_i\in X$ (the degree $n$ component of $T(X)$ is spanned by words of length $n$); $\shuffle$ is the shuffle product recursively defined by
$$y_1\dots y_n\shuffle z_1\dots z_m:=y_1\cdot (y_2\dots y_n\shuffle z_1\dots z_m)+z_1\cdot (y_1\dots y_n\shuffle z_2\dots z_m),$$
where $\cdot$ stands for the concatenation of words, $y_1\dots y_n\cdot z_1\dots z_m:=y_1\dots y_nz_1\dots z_m$; $\Delta$ is the deconcatenation coproduct defined by
$$\Delta(y_1\dots y_n):=\sum\limits_{i=0}^ny_1\dots y_i\otimes y_{i+1}\dots y_n.$$

The symmetric group $S_n$ acts on $T(X)$ as the null map on words of length $m\not=n$ and by permuting the letters of the words of length $n$:
$$\sigma(y_1\dots y_n):=y_{\sigma^{-1}(1)}\dots y_{\sigma^{-1}(n)}.$$
One can show that $\bigoplus\limits_{n=0}^\infty\Q[S_n]$ is a subalgebra of the algebra of linear endomorphisms of $T(X)$ equipped with the convolution product $\ast$. The subalgebra $\desc$ generated as an algebra by the identity maps $1_n$ (seen as words $1_n=12\ldots n$) is called the descent algebra. It is a graded algebra with graded components $\desc_n:=\desc\cap \Q[S_n]$.
It is a free algebra over the elements $1_n$ and as such canonically identifies with the algebra of noncommutative symmetric functions $\NSym$ which is, by definition the free associative graded algebra over a set of graded generators. (In the rest of the article we won't distinguish between $\desc$ and $\NSym$.) It has therefore a basis $D_\mu$ indexed by positive integer sequences  $\mu=(\mu_1,\dots,\mu_k)$, where 
$$D_\mu:=1_{\mu_1}\ast \dots \ast 1_{\mu_k}.$$
Recall that the sequence $\mu$ is also called a composition of the integer $|\mu|=\mu_1+\cdots +\mu_k$.
Writing $S(\mu):=\{\mu_1,\mu_1+\mu_2,\dots,\mu_1+\dots+\mu_{k-1}\}$, one shows (see \cite{Garsia1985}) that $D_\mu=\sum\limits_{\sigma\in D_{\subseteq S(\mu)}}\sigma$. 
Thus, $\desc$ has a basis formed by the $\{D_\mu \ | \  \mu \rm { \ a \ composition }\}$ and, by a standard M\"obius inversion argument, it also has a basis formed by sums of elements in the descent classes in symmetric groups. 

The subspace $\desc_n$ happens to be stable by the product in $\Q[S_n]$ or, equivalently, the composition of linear endomorphisms of $T(X)$. The formula reads, for $\mu=(\mu_1,\dots,\mu_n)$ and $\beta=(\beta_1,\dots,\beta_m)$ two compositions of $n$:
\begin{equation}\label{prodfla}
D_\beta\circ D_\mu:=\sum\limits_{c(\nu)=\mu\atop r(\nu)=\beta}D_{\omega(\nu)},
\end{equation}
where $\nu=(\nu_{i,j})$ runs over $m\times n$ matrices with nonnegative integer coefficients, $c(\nu)$ is the sequence of length $n$ of columns' coefficients sums, $r(\nu)$ is the sequence of length $m$ of lines' coefficients sums. The symbol $\omega(\nu)$ stands for 
$$\omega(\nu)=p(\nu_{1,1},\nu_{1,2},\dots,\nu_{1,n},\dots,\nu_{m,n}),$$
where the map $p$ acts as the packing of nonnegative integer sequences (the removal of their null entries). For example, $p(3,0,6,0,3):=(3,6,3)$.
In words, $\omega(\nu)$ stands for the packing of the word obtained by reading the coefficients of the matrix $\nu$ in the reading order.

The fact that the deconcatenation coproduct obviously satisfies $\sum\limits_{i=0}^n(1_i\otimes 1_{n-i})\circ \Delta=\Delta\circ 1_n$ has important implications that automatically follow from general properties of the behaviour of bialgebra endomorphisms: see \cite[Sect. 3.3]{Cartier2021}. First, as $\desc$ is a free associative algebra over the $1_n$, $n\geq 1$, the coproduct $\blacktriangle$, defined on generators by
$$\blacktriangle(1_n):=\sum\limits_{i=0}^n(1_i\otimes 1_{n-i}),$$ defines a coalgebra structure on all of $\desc$ by
\begin{equation}\label{eqq3}
\blacktriangle (a\ast b)=\blacktriangle(a)\ast\blacktriangle(b).
\end{equation}
The following identities also automatically hold for arbitrary elements $a,b,c$ in $\desc=\NSym$:
\begin{equation}\label{eqq1}
(a\ast b)\circ c=(a\circ c^{(1)})\ast(b\circ c^{(2)}),\end{equation}
\begin{equation}\label{eqq2}
\blacktriangle (a\circ b)=\blacktriangle(a)\circ\blacktriangle(b).\end{equation}
In words: $(\desc,\ast,\blacktriangle)$ is a graded connected cocommutative Hopf algebra; $(\desc,\circ,\blacktriangle)$ is a cocommutative bialgebra (without a unit, as its natural unit would be the infinite sum $\sum\limits_{n=0}^\infty 1_n$ that does not belong to $\desc$ but only to its completion $\prod\limits_{n=0}^\infty \desc_n$); the product $\circ$ is distributive over $\ast$ (in the bialgebraic sense provided by Eq. (\ref{eqq1})).

The scope of the present article will be to extend these constructions from the linear span of compositions of integers  to a two-dimensional analogue, using the product formula Eq. (\ref{prodfla}), the bialgebra structure associated to the composition product implied by Eq. (\ref{eqq2}), and the other properties of $\desc$ as a guideline.

Before coming to these developments we briefly recall how the previous results dualise, following \cite{Malvenuto1995} and \cite{Foissy40}.

The dual of $\desc$ (in the graded sense), $\bigoplus\limits_{n\in\N}\desc_n^\ast$ is the graded vector space of quasi-symmetric functions $\QSym$ (see \cite{Malvenuto1995}).
Let us define the dual basis of the basis of the $D_\mu$ by introducing the pairing $\langle D_\mu,\nu\rangle:=\delta_{\mu, \nu}$, where $\mu$ and $\nu$ run over compositions.
As a vector space, $\QSym$ is thus generated by compositions $\nu=(\nu_1,\ldots,\nu_k)$ .

The product dual to $\blacktriangle$ is called the quasi-shuffle product, it is written $\squplus$ and is explicitly given by
\begin{align*}
(\nu_1,\ldots,\nu_k)\squplus (\nu_{k+1},\ldots,\nu_{k+l})&=\sum_{\sigma \in \qsh(k,l)}
\left(\sum_{\sigma(i)=1} \nu_i,\ldots, \sum_{\sigma(i)=\max(\sigma)}\nu_i\right),
\end{align*}
where $\qsh(k,l)$ is defined as follows.

\begin{defi}
Let $k,l\in \N$. A $(k,l)$-quasi-shuffle is a surjective map $\sigma:[k+l] \longrightarrow [n]$ such that
\begin{align*}
\sigma(1)<\ldots<\sigma(k),&&\sigma(k+1)<\ldots<\sigma(k+l).
\end{align*}
The set of $(k,l)$-quasi-shuffles is denoted $\qsh(k,l)$.
\end{defi}

\begin{example}\begin{align*}
\qsh(1,1)&=\{(11),(12),(21)\}\\
\qsh(2,1)&=\{(121),(122),(123),(132),(231)\}\\
\qsh(1,2)&=\{(112),(123),(212),(213),(312)\}\\
\qsh(2,2)&=\left\{\begin{array}{c}
(1212),(1213),(1223),(1234),(1312),(1323),\\
(1324),(1423),(2312),(2313),(2314),(2413),(3412)
\end{array}\right\}.
\end{align*}\end{example}

It is common to give a polynomial representation of $\QSym$ as a subalgebra of the algebra of formal power series in a countable set of variables $X=\{x_1,\dots,x_n,\dots \}$ whose linear basis are the (generalized) polynomials $M_\nu:=\sum\limits_{i_1<\dots <i_k}x_{i_1}^{\nu_1}\dots x_{i_k}^{\nu_k}$,
where $\nu=(\nu_1,\ldots,\nu_k)$ runs over compositions, the so--called monomial quasi-symmetric functions in the sense of Gessel \cite{Gessel1984}. Using the extended notation $M_\gamma:=\sum\limits_{i=1}^n\lambda_i M_{\mu_i}$ when 
$\gamma:=\sum\limits_{i=1}^n\lambda_i \mu_i$ is a linear combination of compositions $\mu_i$,
one obtains the product formula $M_\nu M_\mu=M_{\nu\squplus\mu}$, where $\nu$ and $\mu$ are two compositions.
We won't distinguish later in the article between the presentation of $\QSym$ as generated by compositions or using its polynomial presentation, excepted for the use of one notation or the other.

As $\desc=\NSym$ is free associative, the coproduct in $\QSym$ dual to the product in $\desc$ is the deconcatenation coproduct
\begin{align*}
\Delta(a_1,\ldots,a_n)&=\sum_{k=0}^n (a_1,\ldots,a_k)\otimes (a_{k+1},\ldots,a_n).
\end{align*}
Thus, $(\QSym,\squplus,\Delta)$ is a Hopf algebra, the Hopf algebra of quasi-symmetric functions, dual to $(\desc=\NSym,\ast,\blacktriangle)$. 
As $(\desc_n,\circ)$ is a finite dimensional algebra, the bialgebra structure (without unit) $(\desc,\circ,\blacktriangle)$ also dualizes to a bialgebra structure (with a unit and a counit) $(\QSym,\squplus,\delta)$, where we write $\delta$ for the coproduct dual to the composition product in $\desc$.
The same formulas and ideas apply to the polynomial realisation  (however in that case one can alternatively use the notions of products of alphabets to construct $\Delta$ and $\delta$, as in \cite{Malvenuto1995}). 

\subsection{Double bialgebras}

\label{sectiondouble}
The full structure of $\QSym$ is best encoded in the notion of double bialgebra \cite{Foissy40} that will be essential later on in the present article. As all the double bialgebras we will consider will be commutative as algebras, we require commutativity of the product in their definition although it is not required in full generality.

\begin{defi}
A (commutative) double bialgebra is a family $(A,m,\Delta,\delta)$ such that:
\begin{itemize}
\item $(A,m,\delta)$ is a commutative bialgebra.
\item $(A,m,\Delta)$ is a commutative bialgebra in the category of right comodules over $(A,m,\delta)$, with the coaction given by $\delta$
itself.
\end{itemize}
In particular, $\delta$ and $\Delta$ satisfies the following compatibility:
\[(\Delta \otimes \id)\circ \delta=(\id \otimes \id \otimes m)\circ (\id \otimes c\otimes \id)\circ (\delta \otimes \delta) \circ \Delta,\]
where $c:A\otimes A\longrightarrow A\otimes A, c(a\otimes b):=b\otimes a$ is the usual flip. The counit $\varepsilon_\Delta$ of $\Delta$ satisfies the following compatibility with $\delta$: for any $x\in A$,
\[(\varepsilon_\Delta \otimes \id)\circ \delta(x)=\varepsilon_\Delta(x)1_A.\] 
\end{defi}
\begin{example}
A simple example of such an object is $\K[X]$, with its usual algebra structure and its two multiplicative coproducts defined by
\[\Delta(X)=X\otimes 1+1\otimes X,\hspace{1cm}\delta(X)=X\otimes X.\]
\end{example}
\begin{example}
Another important example is given by $\QSym$, with the quasi-shuffle product $\squplus$ and its two coproducts $\Delta$ and $\delta$, as it follows from Equations (\ref{eqq1},\ref{eqq2},\ref{eqq3}).
\end{example}
\begin{theo}
\label{theoinvariant}
Let $(A,m,\Delta,\delta)$ be a double bialgebra, such that $(A,m,\Delta)$ is connected.
Then there exists a unique double bialgebra morphism $\phi_A:(A,m,\Delta,\delta)\longrightarrow (\K[X],m,\Delta,\delta)$.
\end{theo}
\begin{example}
The unique double bialgebra morphism $\phi_\QSym$ from $(\QSym,\squplus,\Delta,\delta)$ to $(\K[X],m,\Delta,\delta)$ is given by
\[\phi_\QSym((a_1\ldots a_n))=H_n(X)=\frac{X(X-1)\ldots (X-n+1)}{n!},\]
for any composition $(a_1\ldots a_n)$.
\end{example}
%
We refer to \cite{Foissy40} for details and proofs.

\subsection{Truncation and shuffles}
We briefly discuss in this section two variants of the main constructions presented above that will be generalised later on to the matrix case. The first is a truncation governed by the interpretation of Solomon's descent algebras as a noncommutative version of the representation rings of the symmetric groups, the second is a base change that is a particular case of an isomorphism due to Hoffman.

Let us start with the truncation. In the theory of symmetric functions, one usually considers symmetric functions over a countably infinite alphabet $X$.  The main theorem of the theory states that the bicommutative graded connected Hopf algebra of symmetric functions is a free commutative algebra generated by elementary symmetric functions (or complete symmetric, or power sums...). When the alphabet is finite of cardinality $n$, the same result holds but with a finite family (of cardinality $n$) of generators: the same than in the infinite dimension case, but with a generating family restricted to polynomials of degree less or equal to $n$. If one adopts instead the point of view of
symmetric group representations, the bicommutative graded connected Hopf algebra defined as the direct sum of the representation rings of the various symmetric groups together with the iterated product and coproduct obtained by induction and restriction from Young subgroups, is isomorphic to the Hopf algebra of symmetric functions. Here, the truncation amounts to considering only representations of $S_k$, for $k>n$, induced from representations of Young subgroups $S_{p_1}\times\dots\times S_{p_l}$ with $p_1+\dots+p_l=k$ and the $p_i$ less or equal to $n$. The same results would hold for symmetric group characters.

In the noncommutative framework, this leads to:
\begin{defi}
The $n$-truncated Hopf algebra of noncommutative symmetric functions is the sub--Hopf algebra $\NSym_n$ of $\NSym$ freely generated as an associative algebra by the identity permutations $1_k,\ k\leq n$.
\end{defi}
We let the reader check that the definition is consistent, this follows from the fact that coproducts of the $1_k$ read $\blacktriangle(1_k)=\sum\limits_{i=0}^k1_i\otimes 1_{k-i}$.

Another construction is particularly important in the theory of descents. It is fundamental for example for its applications to the theory of free Lie algebras. It reads as follows. The family of the $1_k\in S_n$ is group-like in the completion of $\NSym$ for the topology induced by the grading:
 $$\Delta\left(\sum\limits_{n=0}^\infty1_n\right)=\sum\limits_{n=0}^\infty1_n\otimes \sum\limits_{n=0}^\infty1_n.$$
Setting $I:=\sum\limits_{n=0}^\infty1_n$, it follows that $\log(I)$ is a primitive element (still in the completion of $\NSym$). Its graded components are known as the canonical, Solomon, or Eulerian idempotents and usually denoted $e_n^1$. By a triangularity argument, $\NSym$ is, as a cocommutative Hopf algebra, the free associative algebra over the primitive elements $e_n^1$, or, equivalently, the enveloping algebra of the free Lie algebra over the $e_n^1$.

The dual statement (that is, for $\QSym$) is a particular case of the Hoffman isomorphism (the case $V=\Q \N^\ast$, the group algebra of the positive integers, with $V$ as below). We present the isomorphism \cite{Hoffman2000,Hoffman2020}, as we will use its Corollary \ref{corbattages} later on. 

Let $V$ be a commutative algebra. For any $v,z\in V$, we denote by $\|vz\|$ their product in $V$.  
As above for $T(X)$, let $T(V)$ stand for the space of tensors over the vector space $V$, $\shuffle$  for the shuffle product, and $\Delta$ for the deconcatenation coproduct. The quasi-shuffle product $\squplus$ is inductively defined (using a word notation as for $T(X)$) by
 $$v_1\dots v_n\squplus z_1\dots z_m:=v_1\cdot(v_2\dots v_n\squplus z_1\dots z_m)+z_1\cdot (v_1\dots v_n\squplus z_2\dots z_m)$$
$$+\|v_1z_1\|\cdot(v_2\dots v_n\squplus z_2\dots z_m),$$
where $\cdot$ denotes here the concatenation product in $T(V)$. 

\begin{theo}
The following map is an isomorphism of bialgebras:
\[\Exp:\left\{\begin{array}{rcl}
(T(V),\shuffle,\Delta)&\longrightarrow&(T(V),\squplus,\Delta)\\
w&\longmapsto&\displaystyle \sum_{w=w_1\ldots w_k} \frac{1}{\ell(w_1)!\ldots \ell(w_k)!}\|w_1\| \ldots \|w_k\|,
\end{array}\right.\]
where for any word $u\in T(V)$, $\ell(u)$ is its length and $\|u\|$ the product of its letters in $V$.
\end{theo}

\begin{notation}
For any $k,l\geq 0$, we denote by $\sh(k,l)$ the set of $(k,l)$-shuffles, that is to say of injective $(k,l)$-quasi-shuffles.
For any $k\geq 1$, we denote by $\inc(k)$ the set of increasing surjective maps $\sigma:[k]\twoheadrightarrow [l]$, i.e. 
\[\sigma(1)\leq  \ldots\leq \sigma(k).\]
If $\sigma:[k]\longrightarrow [l]$ is a map, we put 
\[\sigma!=\prod_{i=1}^l |\sigma^{-1}(i)|!.\]
For any $k,l\geq 0$,
\begin{align*}
|\sh(k,l)|&=\frac{(k+l)!}{k!l!},&|\inc(k)|&=2^{k-1}. 
\end{align*}
\end{notation}

\begin{example}
\begin{align*}
\sh(1,1)&=\{(12),(21)\}\\
\sh(2,1)&=\{(123),(132),(231)\}\\
\sh(1,2)&=\{(123),(213),(312)\}\\
\sh(2,2)&=\{(1234),(1324),(1423),(2314),(2413),(3412)\}\\
\inc(1)&=\{(1)\}\\
\inc(2)&=\{(11),(12)\}\\
\inc(3)&=\{((111),(112),(122),(123)\}\\
\inc(4)&=\{(1111),(1112),(1122),(1123),(1222),(1223),(1233),(1234)\}.
\end{align*}
\end{example}
\begin{cor}\label{corbattages}
For any $k,l\geq 0$, in the vector space generated by the maps $\sigma:[k+l]\twoheadrightarrow [p]$ with $p\leq k+l$,
\[\sum_{\substack{\sigma'\in \inc(k)\\ \sigma''\in \inc(l)\\ \tau \in \qsh(\max(\sigma'),\max(\sigma''))}}\frac{\tau\circ(\sigma' \otimes \sigma'')}{\sigma'!\sigma''!}
=\sum_{\substack{\sigma\in \sh(k,l)\\ \tau\in \inc(k+l)}} \frac{\tau\circ \sigma}{\tau!}. \]
\end{cor}

\begin{proof}
We choose $V=\Q[X_1,\ldots,X_{k+l}]$, whose product is denoted here by $\star$, in order not to be confused with the concatenation product of $T(V)$. 
Let us apply Hoffman's isomorphism, which can be rewritten as
\begin{align*}
&\forall v_1,\ldots,v_n \in V,&\Exp(v_1\ldots v_n)&=\sum_{\sigma \in \inc(k)}\frac{1}{\sigma!}\left(\prod_{\sigma(i)=1}^\star v_i\right)
\ldots \left(\prod_{\sigma(i)=\max(\sigma)}^\star v_i\right).
\end{align*}
Then
\begin{align*}
&\Exp(X_1\ldots X_k)\squplus \Exp(X_{k+1}\ldots X_{k+l})\\
&=\sum_{\substack{\sigma'\in \inc(k)\\ \sigma''\in \inc(l)\\ \tau \in \qsh(\max(\sigma'),\max(\sigma''))}}\frac{1}
{\sigma'!\sigma''!}
\left(\prod_{\tau\circ(\sigma' \otimes \sigma'')=1}^\star X_i \right) \ldots \left(\prod_{\tau\circ(\sigma' \otimes \sigma'')=
\max(\tau)}^\star X_i \right),\\
&=\Exp(X_1\ldots X_k\shuffle X_{k+1}\ldots X_{k+l})\\
&=\sum_{\substack{\sigma\in \sh(k,l)\\ \tau\in \inc(k+l)}} \frac{1}{\tau!}
\left(\prod_{\tau \circ \sigma(i)=1}^\star X_i\right)\ldots
\left(\prod_{\tau \circ \sigma(i)=\max(\tau)}^\star X_i\right).
\end{align*}
Noticing that the words in monomials appearing in this equality are linearly independent, we deduce the result. 
\end{proof}

\section{A bialgebra of integer matrices}\label{sec:bim}

\subsection{Construction}

\begin{defi}
\begin{enumerate}
\item We denote by $\calM_{k,l}(\N)$ the set of matrices with natural coefficients with $k$ rows and $l$ columns, and by $\mat$ the set 
\[\mat=\calM_{0,0}(\N) \cup \bigcup_{k,l\geq 1} \calM_{k,l}(\N),\]
where by convention $\calM_{0,0}(\N)$ is reduced to a single element, denoted by $1$. 
\item If $M=(m_{i,j})_{\substack{1\leq i\leq k\\1\leq j\leq l}}\in \mat$, the number $k$ of rows of $M$ is denoted by $\row(M)$ and the number $l$ of columns of $M$ by $\col(M)$. The weight of $M$ is 
\[\omega(M)=\sum_{i=1}^{\row(M)}\sum_{j=1}^{\col(M)} m_{i,j}=\begin{pmatrix}
1&\ldots&1
\end{pmatrix}M\begin{pmatrix}
1\\ \vdots\\1
\end{pmatrix}.\]
By convention, $\row(1)=\col(1)=\omega(1)=0$.
\item The $\Q$-vector space spanned by $\mat$ is denoted by $\bfH_\mat$.
\end{enumerate} \end{defi}

\begin{prop}
Let $M,M' \in \mat$. We define the matrix $M\searrow M'$ as the block matrix
\[M\searrow M'=\begin{pmatrix}
M&0\\
0&M'
\end{pmatrix}.\]
This product is bilinearly extended to $\bfH_\mat$, making it an associative, unitary, noncommutative algebra. Its unit is the empty matrix $1$.
\end{prop}

\begin{proof}
Immediate.
\end{proof}

For any $k,l\geq 1$, the monoid $(\calM_{k,l}(\N),+)$ is locally finite, which means that for any $M\in \calM_{k,l}(\N)$,
the set of pairs $(M',M'')\in \calM_{k,l}(\N)^2$ such that $M=M'+M''$ is finite. 
\begin{defi}\label{deficoproduit}
We define a map $\blacktriangle:\bfH_\mat\longrightarrow \bfH_\mat\otimes \bfH_\mat$ by $\blacktriangle(1)=1$ and, for any $k,l\geq 1$,
\begin{align*}
&\forall M\in \calM_{k,l}(\N),&\blacktriangle(M)&=\sum_{\substack{(M',M'')\in \calM_{k,l}(\N)^2,\\ M=M'+M''}} M'\otimes M''.
\end{align*}
Then $(\bfH_\mat,\searrow,\blacktriangle)$ is a bialgebra.
\end{defi}

\begin{proof}
As for any locally finite monoid, the associativity of $(\calM_{k,l}(\N),+)$ induces the coassociativity of $\blacktriangle$. The counit is given by
\begin{align*}
&\forall M\in \mat,&\varepsilon_\blacktriangle(M)=\begin{cases}
1\mbox{ if }\omega(M)=0,\\
0\mbox{ otherwise}.
\end{cases}
\end{align*}
Let $M_1,M_2\in \mat$, and let $M',M''\in \mat$ such that $M_1\searrow M_2=M'+M''$. According to the form of $M_1\searrow M_2$ as a block-diagonal matrix, there exist $M'_1$, $M'_2$, $M''_1$, $M''_2$ in $\mat$ such that
$M'=M'_1\searrow M'_2$, $M''=M''_1\searrow M''_2$ and $M_1=M'_1+M''_1$, $M_2=M'_2+M''_2$. Hence,
\[\blacktriangle(M_1\searrow M_2)=\sum_{\substack{M_1=M'_1+M''_1,\\ M_2=M'_2+M''_2}}M'_1\searrow M'_2 \otimes M''_1\searrow M''_2=\blacktriangle(M_1)\searrow \blacktriangle(M_2).\]
So $(\bfH_\mat,\searrow,\blacktriangle)$ is a bialgebra.
\end{proof}

\begin{remark}
\begin{enumerate}
\item The bialgebra $(\bfH_\mat,\searrow,\blacktriangle)$ is graded by the weight: 
\begin{itemize}
\item If $M,M'\in \mat$, then $\omega(M\searrow M')=\omega(M)+\omega(M')$.
\item if $M',M''\in \mat$, then $\omega(M'+M'')=\omega(M')+\omega(M'')$. 
\end{itemize}
However, the homogeneous components of $\bfH_\mat$ are not finite-dimensional. 
\item The bialgebra $(\bfH_\mat,\searrow,\blacktriangle)$ is not a Hopf algebra: for any $k,l\geq 1$, the zero matrix $0_{k,l}$ is a group-like element, with no inverse. 
\end{enumerate}
\end{remark}

\begin{prop}\label{proptransposition}
The following map defined on the linear basis of $\bfH_\mat$ by the usual matrix transposition is an involutive bialgebra automorphism of $(\bfH_\mat,\searrow,\blacktriangle)$:
\[T:\left\{\begin{array}{rcl}
\bfH_\mat&\longrightarrow&\bfH_\mat\\
M\in \mat&\longmapsto&M^\top.
\end{array}\right.\]
\end{prop}

\begin{proof}
For any $M,M'\in \mat$, $(M\searrow M')^\top=M^\top \searrow M'^\top$: $T$ is an algebra morphism.
If $M,M',M''\in \calM_{k,l}(\N)$, then $M=M'+M''$ if, and only if, $M^\top=M'^\top+M''^\top$: this implies that $T$ is a coalgebra morphism. Moreover, $T$ is clearly involutive. 
\end{proof}

\subsection{Quotient to packed matrices}

\begin{defi}
A packed matrix is an element $M=(m_{i,j})_{\substack{1\leq i\leq k\\1\leq j\leq l}}\in \calM_{k,l}(\N)$ such that $M$ has no zero row and no zero column.
By convention, the empty matrix $1$ is packed.
The set of packed matrices is denoted by $\pack$ and the set of packed matrices of weight $n$ is denoted by $\pack_n$ for any $n\in \N$.
\end{defi}

The number of packed matrices of a given weight is given by Entry A120733 of the OEIS \cite{Sloane}.
\begin{align*}
\begin{array}{|c||c|c|c|c|c|c|c|c|c|c|c|}
\hline n&0&1&2&3&4&5&6&7&8&9&10\\
\hline\hline |\pack_n|&1&1&5&33&281&2961&37277&546193&9132865&171634161&3581539973\\
\hline\end{array}\end{align*}

\begin{defi}
For any $M\in \mat$, we denote by $p(M)$ the matrix obtained by the deletion of the zero rows and columns of $M$. This defines a map $p:\mat\longrightarrow \pack$.
Moreover, $p\circ p=p$ and $\omega \circ p=\omega$. 
\end{defi}

\begin{remark}
In particular, for any $k,l\geq 1$, $p(0_{k,l})$ is the empty matrix $1$. 
\end{remark}

\begin{prop}
We put $I=Span(M-p(M)\mid M\in \mat)$. Then $I$ is a biideal of $(\bfH_\mat,\searrow,\blacktriangle)$. 
The quotient is denoted by $(\bfH_\pack,\searrow,\blacktriangle)$. Moreover, $(\bfH_\pack,\searrow,\blacktriangle)$ is a Hopf algebra, graded by the weight.
\end{prop}

\begin{proof}
Let $M,N\in \mat$. Then
\begin{align*}
p(M\searrow N)&=p(p(M)\searrow N)=p(M\searrow p(N))=p(M)\searrow p(N).
\end{align*}
Therefore,
\[(M-p(M))\searrow N=\underbrace{(M\searrow N-p(M)\searrow p(N))}_{=M\searrow N-p(M\searrow N)}
-(\underbrace{p(M)\searrow N-p(M)\searrow p(N)}_{=p(M)\searrow N-p(p(M)\searrow N)})\in I,\]
so $I$ is a right ideal for $\searrow$. Similarly, $I$ is a left ideal.\\

Let $M,M',M''\in \mat$, such that $M'+M''=M$. We denote by $q_M(M')$ and  by $q_M(M'')$ the matrices obtained by deleting in $M'$ and $M''$ the rows and the columns which are deleted in $M$ to obtained $p(M)$ (these rows are also zero in $M$ and $M'$). Then
\begin{align*}
\blacktriangle(M-p(M))&=\sum_{M=M'+M''}M'\otimes M''-\sum_{p(M)=N'+N''} N'\otimes N''\\
&=\sum_{M=M'+M''} M'\otimes M''-q_M(M')\otimes q_M(M'')\\
&=\sum_{M=M'+M''}(M'-q_M(M'))\otimes M''+q_M(M')\otimes (M''-q_M(M'')).
\end{align*}
Moreover, $p(M')=p(q_M(M'))$, so
\[M'-q_M(M')=(M'-p(M'))-(q_M(M')-p(q_M(M')))\in I.\]
Similarly, $M''-q_M(M'')\in I$: we proved that $\blacktriangle(I)\subseteq I\otimes \bfH_\mat\otimes \bfH_\mat \otimes I$,
so $I$ is a biideal. \\

As for the counit, for any $M\in \mat$, $\varepsilon_\blacktriangle(M)=\varepsilon_\blacktriangle(p(M))$, so $I\subseteq \ker(\varepsilon_\blacktriangle)$.\\

For any $M\in \mat$, $\omega(M)=\omega(p(M))$, so $I$ is a graded biideal, and consequently $\bfH_\mat/I$ is a graded bialgebra. Its homogeneous component of weight 0 is reduced to $\Q 1$, so the bialgebra $\bfH_\mat/I$ is connected:
therefore, it is a Hopf algebra. 
\end{proof}

Let us give a description of $\bfH_\pack$.

\begin{prop}
The bialgebra $(\bfH_\pack, \searrow,\blacktriangle)$ has for basis the set $\pack$ of matrices without zero rows or columns. The product is given by
\begin{align*}
&\forall M,M'\in \pack,&M\searrow M'=\begin{pmatrix}
M&0\\
0&M'
\end{pmatrix}.
\end{align*}
The coproduct is given by
\begin{align*}
&\forall M\in \pack,&\blacktriangle(M)&=\sum_{\substack{M,M'\in \mat,\\ M=M'+M''}}p(M')\otimes p(M'').
\end{align*}
\end{prop}

\begin{proof}
As $p$ (linearly extended to $\bfH_\mat$) is a projection on the space spanned by packed matrices, we obtain
\[\bfH_\mat=\Q \pack\oplus I.\]
Therefore, as $\Q$-vector spaces, $\bfH_\pack=\bfH_\mat/I$ and $\Q\pack$ can be identified. So $\bfH_\pack$
has for basis the set of (classes modulo $I$ of) packed matrices. Moreover, for any $M,M'\in \mat$,
\[M\searrow M'\in \pack\Longleftrightarrow M,M'\in \pack.\]
In particular, $\Q\pack$ is a subalgebra of $(\bfH_\mat,\searrow)$, so the identification of $\bfH_\pack$ and $\Q\pack$
 is in fact an algebra isomorphism, which gives the formula for the product of $\bfH_\pack$.  Finally, if $M\in \pack$, in $\bfH_\pack$ we have
\begin{align*}
\blacktriangle(M)&=(p\otimes p)\circ \blacktriangle(M)=(p\otimes p)\left(\sum_{\substack{M,M'\in \mat,\\ M=M'+M''}}M'\otimes M''\right). \qedhere
\end{align*}
\end{proof}

\begin{example} Let $a,b\geq 1$.
\begin{align*}
\blacktriangle((a))&=(a)\otimes 1+1\otimes (a)+\sum_{a'=1}^{a-1}(a')\otimes (a-a'),\\
\blacktriangle\left(\begin{pmatrix}
a\\ b
\end{pmatrix}\right)&=\begin{pmatrix}
a\\ b
\end{pmatrix}\otimes 1+1\otimes \begin{pmatrix}
a\\ b
\end{pmatrix}+\sum_{a'=1}^{a-1}\begin{pmatrix}
a'\\b
\end{pmatrix}\otimes (a-a')+\sum_{a'=1}^{a-1}(a') \otimes \begin{pmatrix}
a-a'\\ b
\end{pmatrix}\\
&+\sum_{b'=1}^{b-1}\begin{pmatrix}
a\\b'
\end{pmatrix}\otimes (b-b')+\sum_{b'=1}^{b-1}(b')\otimes \begin{pmatrix}
a\\b-b'
\end{pmatrix}\\
&+\sum_{a'=1}^{a-1}\sum_{b'=1}^{b-1}\begin{pmatrix}
a'\\b'
\end{pmatrix}\otimes \begin{pmatrix}
a-a'\\b-b'
\end{pmatrix}+(a)\otimes (b)+(b)\otimes (a). 
\end{align*}
\end{example}

The bialgebra $(\bfH_\pack,\searrow,\blacktriangle)$ is graded by the weight. As the unique packed matrix of weight $0$ is the empty matrix 1, this bialgebra is connected: consequently, it is a Hopf algebra. We can compute its antipode by Takeuchi's formula \cite{Takeuchi1971}. 
This gives:

\begin{prop}
The antipode of the Hopf algebra $(\bfH_\pack,\searrow,\blacktriangle)$ is given by
\begin{align*}
&\forall M\in \pack_n,&S(M)=\sum_{k=1}^n (-1)^k\sum_{\substack{M=M_1+\cdots+M_k,\\ \forall i\in [k],\: \omega(M_i)\neq 0}} p(M_1)\searrow \ldots \searrow p(M_k).
\end{align*}
\end{prop}

\begin{example}
For any $n\geq 1$,
\[S((n))=\sum_{k=1}^n (-1)^k \sum_{\substack{n_1+\cdots+n_k=n,\\ \forall i\in [k],\: n_i>0}}
\begin{pmatrix}
n_1&0&\ldots&0\\
0&\ddots&\ddots&\vdots\\
\vdots&\ddots&\ddots&0\\
0&\ldots&0&n_k
\end{pmatrix}.\]
\end{example}

\begin{prop}
The map $T$ of Proposition \ref{proptransposition} induces an involutive bialgebra morphism of $(\bfH_\pack,\searrow,\blacktriangle)$:
\[T:\left\{\begin{array}{rcl}
\bfH_\pack&\longrightarrow&\bfH_\pack\\
M\in \mat&\longmapsto&M^\top.
\end{array}\right.\]
\end{prop}

\begin{proof}
This comes from the fact that for any $M\in \mat$, $p(M)^\top=p(M^\top)$. 
\end{proof}

\subsection{A bialgebra of row matrices and $\NSym$}

\begin{defi}
We denote by $\mat_1$ the set of row matrices which coefficients are natural integers:
\[\mat_1=\calM_{0,0}(\N)\cup \bigcup_{l=1}\calM_{1,l}(\N).\]
We denote by $\bfH_{\mat_1}$ the space generated by $\mat_1$. It is given a bialgebra structure by the following:
\begin{align*}
&\forall M,M'\in \mat_1,&M\rightarrow M'&=(M M'),\\
&\forall M\in \calM_{1,l}(\N),&\blacktriangle(M)&=\sum_{\substack{M',M''\in \calM_{1,l}(\N),\\ M=M'+M''}}M'\otimes M''.
\end{align*}
The unit is the empty matrix 1. The counit is given by
\begin{align*}
&\forall M\in \mat_1,&\varepsilon_\blacktriangle(M)=\begin{cases}
1\mbox{ if }\omega(M)=0,\\
0\mbox{ otherwise}.
\end{cases}
\end{align*}
\end{defi}

\begin{proof}
Similar to the proof of Definition \ref{deficoproduit}.
\end{proof}

\begin{remark}
By definitition of the product, $\bfH_{\mat_1}$ is the free algebra generated by the elements $(n)$, with $n\geq 0$. The coproduct of such an element is given by
\[\blacktriangle((n))=\sum_{k=0}^n (k)\otimes (n-k).\]
\end{remark}

\begin{prop}
We denote by $\pack_1$ the set of packed row matrices which coefficients are natural integers:
\[\pack_1=\pack_{0,0}\cup \bigcup_{l=1}\pack_{1,l}.\]
We denote by $\bfH_{\pack_1}$ the space generated by $\pack_1$. The packing map $p$ induces a bialgebra structure on $\bfH_{\pack_1}$, such that
\begin{align*}
&\forall M,M'\in \pack_1,&M\rightarrow M'&=(M M'),\\
&\forall M\in \pack_{1,l}(\N),&\blacktriangle(M)&=\sum_{\substack{M',M''\in \calM_{1,l}(\N),\\ M=M'+M''}}p(M')\otimes p(M'').
\end{align*}
The unit is the empty matrix 1. The counit is given by
\begin{align*}
&\forall M\in \pack_1,&\varepsilon_\blacktriangle(M)=\begin{cases}
1\mbox{ if }\omega(M)=0,\\
0\mbox{ otherwise}.
\end{cases}
\end{align*}
 \end{prop}

\begin{proof}
We consider the space $I$ of $\bfH_{\mat_1}$ generated by the elements $M-p(M)$, with $M\in \pack_1$.
It is the ideal of $\bfH_{\mat_1}$ generated by $(0)-1$. Moreover, 
\[\blacktriangle((0)-1)=(0)\otimes (0)-1\otimes 1=((0)-1)\otimes (0)+1\otimes ((0)-1),\]
so it is a biideal of $(\bfH_{\mat_1},\rightarrow,\blacktriangle)$. Identifying $\bfH_{\mat_1}/I$ and $\bfH_{\pack_1}$,
we obtain the announced bialgebra structure on $\bfH_{\pack_1}$. 
\end{proof}

\begin{remark}
As an algebra, $\bfH_{\pack_1}$ is freely generated by the elements $(n)$, with $n\geq 1$. The coproduct of such an element is given by
\[\blacktriangle((n))=1\otimes(n)+\sum_{k=1}^{n-1} (k)\otimes (n-k)+(n)\otimes 1.\]
Therefore,  $(\bfH_{\pack_1},\rightarrow,\blacktriangle)$ is trivially isomorphic to the Hopf algebra $(\desc,\ast,\Delta)=(\NSym,\ast,\Delta)$ of noncommutative symmetric functions, mapping $(n)$ to $1_n$.

As announced previously, the following proposition shows that the map from packed matrices to elements in the descent algebra, which is involved in the description of products of (sum of elements in the) descent classes in Solomon's descent algebras $\desc_n$ can be lifted to a bialgebra map. 
Corollary \ref{corKxy} further below proves a similar statement in the other direction (from descent classes to packed matrices).
As explained in the Introduction, we expect such phenomena to be meaningful for the noncommutative representation theory of the symmetric groups. We have already obtained some results in that direction and intend to come back on this Remark in a forthcoming work.
\end{remark}

\begin{prop}\label{propTheta}
Let $M=(m_{i,j})_{\substack{1\leq i\leq k,\\ 1\leq j\leq l}}\in\mat$. If we delete the zeros from the sequence 
$(m_{1,1},\ldots,m_{1,l},\ldots,m_{k,1},\ldots, m_{k,l})$ obtained reading the rows of $M$, we get a packed row matrix which we denote by $\comp(M)$. The following map is a bialgebra map  from $(\bfH_\pack,\searrow,\blacktriangle)$ to $(\bfH_{\pack_1},\rightarrow,\blacktriangle)$:
\[\Theta:\left\{\begin{array}{rcl}
\bfH_\pack&\longrightarrow&\bfH_{\pack_1}\\
M\in\pack&\longrightarrow&\comp(M).
\end{array}\right.\]
\end{prop}

\begin{proof}
If $M,M'\in \pack$,
\[\Theta(M\searrow M')=\comp\left(\begin{pmatrix}
M&0\\0&M'
\end{pmatrix}\right)=(\comp(M)\comp(M'))=\Theta(M)\rightarrow \Theta(M'),\]
so $\Theta$ is an algebra morphism. Let $M\in \pack$. 
\begin{align*}
(\Theta \otimes \Theta)\circ \blacktriangle(M)&=\sum_{M=M'+M''}\comp \circ p(M')\otimes \comp \circ p(M'')\\
&=\sum_{M=M'+M''}\comp(M')\otimes \comp(M'')\\
&=\sum_{\comp(M)=N'+N''}p(N')\otimes p(N'')\\
&=\blacktriangle\circ \Theta(M). 
\end{align*}
So $\Theta$ is a bialgebra morphism. 
\end{proof}

\begin{prop}
Let us fix $k,l\geq 1$. For any $n\in \N$, we consider the element
\[K'_{k,l}(n)=\sum_{M\in \calM_{k,l}(\N),\: \omega(M)=n} M.\]
Then the following is a bialgebra morphism from $(\bfH_{\mat_1},\rightarrow,\blacktriangle)$ to $(\bfH_\mat,\searrow,\blacktriangle)$:
\[K'_{k,l}:\left\{\begin{array}{rcl}
\bfH_{\mat_1}&\longrightarrow&\bfH_\mat\\
(a_1,\ldots,a_n)&\longmapsto&K'_{k,l}(a_1)\searrow \ldots \searrow K'_{k,l}(a_n).
\end{array}\right.\]
\end{prop}

\begin{proof}
First, note that, as $k$ and $l$ are fixed, the sum defining $K'_{k,l}(n)$ is finite. 
As the algebra $(\bfH_{\mat_1},\rightarrow)$ is freely generated by the elements $(n)$, with $n\geq 0$, $K'_{k,l}$ is the algebra morphism sending $(n)$ to $K'_{k,l}(n)$ for any $n\geq 0$, so it is indeed an algebra morphism.
In order to prove it is a bialgebra morphism, it is enough to prove that for any $n\geq 0$,
\[(K'_{k,l}\otimes K'_{k,l})\circ \blacktriangle((n))=\blacktriangle\circ K'_{k,l}((n)).\]
Indeed,
\begin{align*}
(K'_{k,l}\otimes K'_{k,l})\circ \blacktriangle((n))&=\sum_{i=0}^n K'_{k,l}(i)\otimes K'_{k,l}(n-i)\\
&=\sum_{i=0}^n \sum_{\substack{M',M''\in \calM_{k,l}(\N),\\ \omega(M')=i,\:\omega(M'')=n-i}}M'\otimes M''\\
&=\sum_{\substack{M',M''\in \calM_{k,l}(\N),\\ \omega(M'+M'')=n}}M'\otimes M''\\
&=\sum_{\substack{M\in \calM_{k,l}(\N),\\ \omega(M)=n}}\sum_{\substack{M',M''\in \calM_{k,l}(\N),\\ M=M'+M''}}M'\otimes M''\\
&=\sum_{\substack{M\in \calM_{k,l}(\N),\\ \omega(M)=n}} \blacktriangle(M)\\
&=\blacktriangle\circ K'_{k,l}((n)). \qedhere
\end{align*}
\end{proof}

\begin{cor}\label{corKxy}
Let $x,y\in \Q$. The following defines a bialgebra morphism from $(\bfH_{\pack_1},\rightarrow,\blacktriangle)$
to $(\bfH_\pack,\searrow,\blacktriangle)$:
\[K_{x,y}:\left\{\begin{array}{rcl}
\bfH_{\pack_1}&\longrightarrow&\bfH_\pack\\
(a_1 \ldots a_n)&\longmapsto&\displaystyle \sum_{M_i\in \pack,\: \omega(M_i)=a_i}\left(\prod_{i=1}^n H_{\row(M_i)}(x)
H_{\col(M_i)}(y) \right) M_1\searrow \ldots \searrow M_n.
\end{array}\right.\]
\end{cor}

\begin{proof}
We  start with $(x,y)=(k,l)$ a pair of positive integers. Note that $K'_{k,l}((0))=0_{k,l}$.
Therefore, $(0)-1\in \ker(p\circ K'_{k,l})$, so $K'_{k,l}$ induces a bialgebra map $K_{k,l}$ from the quotient $\bfH_{\pack_1}$
of $\bfH_{\mat_1}$ by $(0)-1$ to $\bfH_\pack$. For any $n\geq 1$,
\begin{align*}
K_{k,l}((n))&=\sum_{M\in \calM_{k,l}(\N),\:\omega(M)=n} p(M)\\
&=\sum_{M\in \pack,\: \omega(M)=n} \binom{k}{\row(M)}\binom{l}{\col(M)} M\\
&=\sum_{M\in \pack,\: \omega(M)=n}H_{\row(M)}(k) H_{\col(M)}(l) M,
\end{align*}
where the two binomial coefficients correspond to the "unpacking" of $M$ as a $k\times l$ matrix by adding zero rows and columns.
By multiplicativity, we obtain the formula for $K_{k,l}(a_1\ldots a_n)$ when $n\geq 2$.\\

For any $x,y\in \Q$, $K_{x,y}$ can be defined as the unique algebra morphism from $\bfH_{\pack_1}$ to $\bfH_\pack$ which sends $(n)$ to 
\[K_{x,y}((n))=\sum_{M\in \pack,\: \omega(M)=n} H_{\row(M)}(x)H_{col(M)}(y) M,\]
so it is indeed an algebra morphism. Let us now fix an element $M\in \bfH_{\pack_1}$. We consider the set
\[E_M=\{(x,y)\in \Q^2\mid (\blacktriangle\circ K_{x,y}-(K_{x,y}\otimes K_{x,y})\circ\blacktriangle)(M)=0\}.\]
This set is determined by a polynomial equation $P_M(X,Y)\in (\bfH_\pack\otimes \bfH_\pack)[X,Y]$.
The first part of this proof shows that $P_M(X,Y)$ vanishes on $(\N\setminus\{0\})^2$, so $P_M(X,Y)=0$.
Therefore, $E_M=\Q^2$: for any $(x,y)\in \Q^2$, $K_{x,y}$ is a bialgebra morphism.
\end{proof}

\begin{example}
\begin{align*}
K_{x,y}((1))&=xy(1),\\
K_{x,y}((2))&=xy(2)+\frac{xy(y-1)}{2}\begin{pmatrix}
1&1
\end{pmatrix}+\frac{x(x-1)y}{2}\begin{pmatrix}
1\\1
\end{pmatrix}\\
&+\frac{x(x-1)y(y-1)}{4}\begin{pmatrix}
1&0\\0&1
\end{pmatrix}+\frac{x(x-1)y(y-1)}{4}\begin{pmatrix}
0&1\\1&0
\end{pmatrix},\\
K_{x,y}(\begin{pmatrix}
1&1
\end{pmatrix})&=x^2y^2\begin{pmatrix}
1&0\\0&1
\end{pmatrix}.
\end{align*}
\end{example}

\begin{remark}
Note that $H_n(1)=1$ if $n=0$ or $1$, and $0$ otherwise. Therefore, for any $(a_1\ldots a_n)\in \pack_1$,
\[K_{1,1}((a_1\ldots a_n))=\begin{pmatrix}
a_1&0&\cdots&0\\
0&\ddots&\ddots&\vdots\\
\vdots&\ddots&\ddots&0\\
0&\ldots&0&a_n
\end{pmatrix}.\]
\end{remark}

\subsection{Truncation and other variants}

Recall that the bialgebra of noncommutative symmetric functions can be truncated: it has subbialgebras freely generated by the $1_k,\ k\leq n$. This corresponds to a natural truncation for symmetric functions and symmetric group representations.

The same process can be used on packed matrices:
\begin{defi}
The bialgebra $(\bfH_{\pack,n},\searrow,\blacktriangle)$ is the subbialgebra of $(\bfH_\pack,\searrow,\blacktriangle)$ spanned linearly by packed matrices whose coefficients are less or equal to $n$.
\end{defi}
The fact that the definition is consistent automatically follows from the definitions of $\searrow$ and $\blacktriangle$ that preserve the property of matrix coefficients being less or equal to $n$.

\begin{remark}
An interesting particular case is the one where all matrix coefficients are equal to 0 or 1, $(\bfH_{\pack,1},\searrow,\blacktriangle)$. In that case, notice that square matrices of size $p\times p$ encode binary relations on $[p]$.
\end{remark}
\begin{remark}
One can also show that permutation matrices form a subbialgebra  of $(\bfH_{\pack,1},\searrow,\blacktriangle)$, as they can be characterized as square matrices with columns and rows with only one nonzero entry equal to 1, a property preserved by $\searrow$ and $\blacktriangle$. This bialgebra is graded connected cocommutative and free associative. In particular it is the enveloping algebra of a free Lie algebra, see \cite{Patras2004}.
\end{remark}

Another interesting phenomenon relates to the graduation. Packed matrices $M$ are canonically bigraded by the pairs $(\row(M), \col(M))\in \N^\ast\times\N^\ast$, and the product $\searrow$ is additive for this bigrading: 
$$(\row(M\searrow M'), \col(M\searrow M'))=(\row(M), \col(M))+(\row(M'), \col(M')).$$
Let us write $\pack_{p,q}$ for the set of $(p,q)$-graded packed matrices.

The coproduct $\blacktriangle$ instead does not preserve the graduation, but it satisfies
$$\blacktriangle (\Q \pack_{p,q}) \subset \bigoplus\limits_{p'+p''\geq p\atop q'+q''\geq q} \Q \pack_{p',q'}\otimes \Q \pack_{p'',q''}.$$
Let us write $\bdeltares$ for its obvious restriction to a map from
$\Q \pack_{p,q}$ to $\bigoplus\limits_{p'+p''= p\atop q'+q''= q} \Q \pack_{p',q'}\otimes \Q \pack_{p'',q''}.$
The new coproduct $\bdeltares$ on $\bfH_{\pack}$ is coassociative and an algebra map (by direct inspection, for degree reasons): it defines a new, bigraded, bialgebra structure on $\bfH_{\pack}$.

\begin{defi}\label{bigrapack}
The bigraded bialgebra of packed matrices is the bialgebra $(\bfH_{\pack},\searrow,\bdeltares)$, with the grading $\gr(M):=(\row(M), \col(M))$.
\end{defi}

\section{Duality}\label{sec:dualty}

\subsection{Graded dual of $\bfH_\pack$}\label{ssect:gdhp}

The bialgebra $(\bfH_\pack,\searrow,\blacktriangle)$ is graded, and its homogeneous components are finite-dimensional. We identify its graded dual with $\bfH_\pack$, through the pairing defined by
\begin{align*}
&\forall M,M'\in \pack,&\langle M,M'\rangle&=\delta_{M,M'}.
\end{align*}
It is a bialgebra, denoted by $(\bfH_\pack,\squplus,\Delta)$.

By duality with the product $\searrow$, we immediately obtain:

\begin{lemma}\label{lemcopr}
For any $M\in \pack$,
\[\Delta(M)=\sum_{\substack{M',M''\in \pack,\\ M=M'\searrow M''}}M'\otimes M''.\]
The counit is given by
\[\varepsilon_\Delta(M)=\delta_{M,1}.\]
\end{lemma}

\begin{notation}
Let $\alpha:[k]\longrightarrow [l]$ be a map. The associated matrix $\mu(\alpha)\in \calM_{l,k}(\N)$ is defined by
\[\mu(\alpha)_{i,j}=\delta_{i,\alpha(j)}.\]
Note that $\mu(\alpha)$ belongs to $\pack$ if and only if $\alpha$ is surjective. Moreover, if $\alpha:[k]\longrightarrow [l]$ and $\beta:[l]\longrightarrow [n]$
are maps, then 
\[\mu(\beta \circ \alpha)=\mu(\beta)\mu(\alpha).\]
\end{notation}

\begin{prop}\label{propqshu}
Let $M,M'\in \pack$, with respectively $k$ and $k'$ rows, and $l$ and $l'$ columns. Then
\[M\squplus M'=\sum_{\sigma \in \qsh(k,k'),\:\tau \in \qsh(l,l')} \mu(\sigma) (M\searrow M')\mu(\tau)^\top.\]
\end{prop}

\begin{proof}
We put
\[M\squplus M'=\sum_{N\in \pack} a_{M,M'}^N N.\]
Then, for any $N\in \pack$,
\begin{align*}
a_{M,M'}^N&=\langle M\squplus M',N\rangle\\
&=\langle M\otimes M',\blacktriangle(N)\rangle\\
&=\sum_{\substack{N', N''\in \mat,\\ N=N'+N''}}\langle M,p(N') \rangle \langle M',p(N'')\rangle.
\end{align*}
Therefore, $a_{M,M'}^N$ is the cardinality of the set
\[A_{M,M'}^N=\{(N',N'')\in \mat\mid N=N'+N'',\: p(N')=M,\:p(N'')=M'\}.\]
We also consider the set
\[B_{M,M'}^N=\{(\sigma,\tau)\in \qsh(k,k')\times \qsh(l,l')\mid \mu(\sigma)(M\searrow M')\mu(\tau)^\top=N\}.\]
If $(N',N'')\in A_{M,M'}^N$, then any of the $k$ first rows of $M\searrow M'$ contributes to a row of $N'$
(and so of $N$) and any of the $k'$ last rows of $M\searrow M'$ contributes to a row of $N''$
(and so of $N$). A similar fact holds for the columns. Therefore, we obtain a pair of quasi-shuffles
$(\sigma,\tau) \in \qsh(k,k')\times \qsh(l,l')$ such that
\begin{align*}
N'&=\mu(\sigma)(M\searrow 0_{k',l'})\mu(\tau)^\top,&N''&=\mu(\sigma)(0_{k,l}\searrow M')\mu(\tau)^\top.
\end{align*} 
Therefore,
\begin{align*}
\mu(\sigma)(M\searrow M')\mu(\tau)^\top=\mu(\sigma)(M\searrow 0_{k',l'})\mu(\tau)^\top+\mu(\sigma)(0_{k,l}\searrow M')\mu(\tau)^\top=N'+N''=N.
\end{align*}
So $(\sigma,\tau)\in B_{M,M'}^N$.
Therefore, this defines a map $f:A_{M,M'}^N\longrightarrow B_{M,M'}^N$. Conversely, if $(\sigma,\tau)\in B_{M,M'}^N$, we put
\begin{align*}
N'&=\mu(\sigma)(M\searrow 0_{k',l'})\mu(\tau)^\top,&N''&=\mu(\sigma)(0_{k,l}\searrow M')\mu(\tau)^\top.
\end{align*} 
Then $p(N')=M$, $p(N'')=M'$ and $N+N'=\mu(\sigma)(M\searrow M')\mu(\tau)^\top=N$: 
this defines a map $g:B_{M,M'}^N\longrightarrow A_{M,M'}^N$. It is immediate to see that $f$ and $g$ are inverse bijections. Hence, 
\begin{align*}
M\squplus M'&=\sum_{N\in \pack} |A_{M,M'}^N|N\\
&=\sum_{N\in \pack} |B_{M,M'}^N|N\\
&=\sum_{\sigma \in \qsh(k,k'),\:\tau \in \qsh(l,l')} \mu(\sigma) (M\searrow M')\mu(\tau)^\top. \qedhere
\end{align*}
\end{proof}

\begin{example}
Let $a,b,c\in \N_{\geq 1}$.
\begin{align*}
\begin{pmatrix}
a
\end{pmatrix}\squplus \begin{pmatrix}
b
\end{pmatrix}&=\begin{pmatrix}
a&0\\0&b
\end{pmatrix}+\begin{pmatrix}
0&b\\a&0
\end{pmatrix}+\begin{pmatrix}
b&0\\0&a
\end{pmatrix}+\begin{pmatrix}
0&a\\b&0
\end{pmatrix}\\
&+\begin{pmatrix}
a&b
\end{pmatrix}+\begin{pmatrix}
a\\b
\end{pmatrix}+\begin{pmatrix}
b&a
\end{pmatrix}+\begin{pmatrix}
b\\a
\end{pmatrix}+\begin{pmatrix}
a+b
\end{pmatrix},\\
\begin{pmatrix}
a
\end{pmatrix}\squplus \begin{pmatrix}
b\\c
\end{pmatrix}&=\begin{pmatrix}
a&0\\0&b\\0&c
\end{pmatrix}+\begin{pmatrix}
0&b\\a&0\\0&c
\end{pmatrix}+\begin{pmatrix}
0&b\\0&c\\a&0
\end{pmatrix}+\begin{pmatrix}
a&b\\0&c
\end{pmatrix}+\begin{pmatrix}
0&b\\a&c
\end{pmatrix}\\
&+\begin{pmatrix}
0&a\\b&0\\c&0
\end{pmatrix}+\begin{pmatrix}
b&0\\0&a\\c&0
\end{pmatrix}+\begin{pmatrix}
b&0\\c&0\\0&a
\end{pmatrix}+\begin{pmatrix}
b&a\\c&0
\end{pmatrix}+\begin{pmatrix}
b&0\\c&a
\end{pmatrix}\\
&+\begin{pmatrix}
a\\b\\c
\end{pmatrix}+\begin{pmatrix}
b\\a\\c
\end{pmatrix}+\begin{pmatrix}
b\\c\\a
\end{pmatrix}+\begin{pmatrix}
a+b\\c
\end{pmatrix}+\begin{pmatrix}
b\\a+c
\end{pmatrix}.
\end{align*}
\end{example}

The quasishuffle product $\squplus$ and the deconcatenation coproduct on $\bfH_\pack$ agree with the ones introduced by J. Diehl and L. Schmitz to define an appropriate notion of discrete surface signatures \cite{Diehl}. We get:
\begin{cor}
By Lemma \ref{lemcopr} and Prop. \ref{propqshu}, $(\bfH_\pack,\squplus,\Delta)$, the graded dual Hopf algebra of $(\bfH_\pack,\searrow,\blacktriangle)$, identifies with the Diehl-Schmitz Hopf algebra of matrix compositions.
\end{cor}

\subsection{A structure result}

\begin{prop}
The triple $(\bfH_\pack,\searrow,\Delta)$ is an infinitesimal bialgebra, in the sense of \cite{Loday2006}.
\end{prop}

\begin{proof}
Let $M,M'\in \pack$. We assume that $M\searrow M'=N\searrow N'$, for two packed matrices $N$ and $N'$. 
If $\row(N)\leq \row(M)$ and $\col(N)>\col(M)$, or if $\row(N)>\row(M)$ and $\col(N)\leq \col(M)$ then $N$ contains a zero column or a zero row, so is not packed: 
this is a contradiction. Therefore, two cases are possible:
\begin{enumerate}
\item $\row(N)\leq \row(M)$ and $\col(N)\leq \col(M)$: then there exists $N''$ such that $M=N\searrow N''$, and $N'=N'' \searrow M'$.
\item $\row(N)\geq \row(M)$ and $\col(N)\geq \col(M)$: then there exists $N''$ such that $M'=N''\searrow N'$ and $N=M\searrow N''$.
\end{enumerate}
These cases are not disjoint, their intersection is the case $\row(N)=\row(M)$ and $\col(N)=\col(M)$, where $M=N$ and $M'=N'$. Therefore,
\begin{align*}
\Delta(M\searrow M')&=\sum_{M=N \searrow N''} N\otimes N''\searrow M'+\sum_{M'=N'' \searrow N'} M\searrow N''\otimes N'-M\otimes M'\\
&=\Delta(M)\searrow (1\otimes M')+(M\otimes 1)\searrow \Delta(M')-M\otimes M'.  \qedhere
\end{align*}
\end{proof}

As $(\bfH_\pack,\searrow,\Delta)$ is connected, by Loday and Ronco's rigidity theorem \cite{Loday2006}:

\begin{cor}
The coalgebra $(\bfH_\pack,\Delta)$ is cofree. The space of its primitive elements has for basis the set of $\searrow$-indecomposable packed matrices,
that is to say nonempty packed matrices $M$ such that
\begin{align*}
&\forall N,N'\in \pack,&M=N\searrow N'\Longrightarrow (N=1 \mbox{ or }N'=1).
\end{align*}
\end{cor}

\begin{remark}
The number of indecomposable packed matrices can be computed with the use of formal series. We obtain
\begin{align*}
\begin{array}{|c||c|c|c|c|c|c|c|c|c|c|}
\hline n&1&2&3&4&5&6&7&8&9&10\\
\hline\hline \dim(\prim(\bfH_\pack)_n)&1&4&24&204&2224&29156&443320&7646684&147367456&3137652676\\
\hline\end{array}\end{align*}
As any commutative connected bialgebra, $(\bfH_\pack,\squplus,\Delta)$ is, as an algebra, isomorphic to a polynomial algebra. It has been studied as such by L. Schmitz and N. Tapia in \cite{Schmitz}, to which we refer.
The number $g_n$ or required generators in each degree $n$ can be computed with the use of formal series. We obtain
\begin{align*}
\begin{array}{|c||c|c|c|c|c|c|c|c|c|c|}
\hline n&1&2&3&4&5&6&7&8&9&10\\
\hline\hline g_n&1&4&28&238&2568&32938& 491700& 8350536& 158944092& 3350761964\\
\hline\end{array}\end{align*}
\end{remark}

\subsection{Dual morphisms}

Recall that $(\QSym,\squplus,\Delta)$, the Hopf algebra of quasi-symmetric functions, is the graded dual of the Hopf algebra $\NSym$, which we have identified with $(\bfH_{\pack_1},\rightarrow,\blacktriangle)$. The Hopf paring between $(\QSym,\squplus,\Delta)$ and $(\bfH_{\pack_1},\rightarrow,\blacktriangle)$ is given by
\begin{align*}
&\forall M,M'\in \pack_1,&\langle M,M'\rangle&=\delta_{M,M'}.
\end{align*}

Firstly, let us transpose the bialgebra morphism $T$  of Proposition \ref{proptransposition}.

\begin{prop}\label{propT2}
The following defines an involutive automorphism of $(\bfH_\mat,\squplus,\Delta)$:
\[T:\left\{\begin{array}{rcl}
\bfH_\mat&\longrightarrow&\bfH_\mat\\
M\in \mat&\longmapsto&M^\top.
\end{array}\right.\]
\end{prop}

Let us now dualize the morphism $\Theta$ of Proposition \ref{propTheta}.

\begin{prop}\label{proptheta2}
The following map is a bialgebra map  from $(\QSym,\squplus,\Delta)$ to $(\bfH_\pack,\squplus,\Delta)$:
\[\theta:\left\{\begin{array}{rcl}
\QSym&\longrightarrow&\bfH_\pack\\
(a_1,\ldots,a_k)&\longmapsto&\displaystyle \sum_{M\in \pack,\: \comp(M)=(a_1,\ldots,a_k)} M.
\end{array}\right.\]
\end{prop}

\begin{example}
Let $a,b,c\in \N_{\geq 1}$.
\begin{align*}
\theta(a)&=\begin{pmatrix}
a
\end{pmatrix},\\
\theta(ab)&=\begin{pmatrix}
a&b
\end{pmatrix}+\begin{pmatrix}
a\\b
\end{pmatrix}+\begin{pmatrix}
a&0\\0&b
\end{pmatrix}+\begin{pmatrix}
0&a\\b&0
\end{pmatrix},\\
\theta(abc)&=\begin{pmatrix}
a&b&c
\end{pmatrix}+\begin{pmatrix}
a&0\\b&c
\end{pmatrix}+\begin{pmatrix}
0&a\\b&c
\end{pmatrix}+\begin{pmatrix}
a&0&0\\
0&b&c
\end{pmatrix}+\begin{pmatrix}
0&a&0\\
b&0&c
\end{pmatrix}+\begin{pmatrix}
0&0&a\\
b&c&0
\end{pmatrix}\\
&+\begin{pmatrix}
a&b\\c&0
\end{pmatrix}+\begin{pmatrix}
a&b\\0&c
\end{pmatrix}+\begin{pmatrix}
0&a&b\\c&0&0
\end{pmatrix}+\begin{pmatrix}
a&0&b\\0&c&0
\end{pmatrix}+\begin{pmatrix}
a&b&0\\0&0&c
\end{pmatrix}\\
&+\begin{pmatrix}
a&0\\0&b\\c&0
\end{pmatrix}+\begin{pmatrix}
a&0\\0&b\\0&c
\end{pmatrix}+\begin{pmatrix}
a&0\\b&0\\0&c
\end{pmatrix}+\begin{pmatrix}
0&a\\b&0\\c&0
\end{pmatrix}+\begin{pmatrix}
0&a\\b&0\\c&0
\end{pmatrix}+\begin{pmatrix}
0&a\\0&b\\c&0
\end{pmatrix}+\begin{pmatrix}
a\\b\\c\end{pmatrix}\\
&+\begin{pmatrix}
a&0&0\\0&b&0\\0&0&c
\end{pmatrix}+\begin{pmatrix}
a&0&0\\0&0&b\\0&c&0
\end{pmatrix}+\begin{pmatrix}
0&a&0\\b&0&0\\0&0&c
\end{pmatrix}+\begin{pmatrix}
0&a&0\\0&0&b\\c&0&0
\end{pmatrix}+\begin{pmatrix}
0&0&a\\b&0&0\\0&c&0
\end{pmatrix}+\begin{pmatrix}
0&0&a\\0&b&0\\c&0&0
\end{pmatrix}.
\end{align*}
\end{example}

\begin{remark}
The number of terms of $\theta(a_1,\ldots,a_n)$ is given by
\[q_n=\sum_{n_1+\cdots+n_k=n} |\qsh(n_1,\ldots,n_k)|,\]
see Entry A101370 of the OEIS \cite{Sloane}.
\begin{align*}
\begin{array}{|c||c|c|c|c|c|c|c|c|c|c|}
\hline n&1&2&3&4&5&6&7&8&9&10\\
\hline\hline q_n &1&4&24&196&2016&24976&361792&5997872&111969552&2324081728\\
\hline \end{array}\end{align*}\end{remark}

Finally, let us dualize the morphisms $K_{x,y}$ of Corollary \ref{corKxy}:

\begin{prop}\label{propkxy2}
Let $x,y\in \Q$. The following defines a bialgebra morphism from $(\bfH_\pack,\squplus,\Delta)$ to $(\QSym,\squplus,\Delta)$ by
\begin{align*}
\kappa_{x,y}&:\left\{\begin{array}{rcl}
\bfH_\pack&\longrightarrow&\QSym\\
M\in \pack&\longmapsto&\displaystyle \sum_{k=1}^\infty\: \sum_{\substack{M=M_1\searrow \ldots \searrow M_k,\\ M_1,\ldots,M_k\neq 1}}\:
\left(\prod_{i=1}^k H_{\row(M_i)}(x)H_{\col(M_i)}(y)\right)(\omega(M_1)\ldots\omega(M_k)).
\end{array}\right.
\end{align*}
It is surjective if, and only if, $x\neq 0$ and $y\neq 0$.
\end{prop}

\begin{example}
Let $a,b\in \N_{\geq 1}$.
\begin{align*}
\kappa_{x,y}\begin{pmatrix}
a
\end{pmatrix}&=xy (a)\\
\kappa_{x,y}\begin{pmatrix}
a&b
\end{pmatrix}&=\frac{xy(y-1)}{2} (a+b)\\
\kappa_{x,y}\begin{pmatrix}
a\\b
\end{pmatrix}&=\frac{x(x-1)y}{2} (a+b)\\
\kappa_{x,y}\begin{pmatrix}
0&a\\b&0
\end{pmatrix}&=\frac{x(x-1)y(y-1)}{4} (a+b)\\
\kappa_{x,y}\begin{pmatrix}
a&0\\0&b
\end{pmatrix}&=\frac{x(x-1)y(y-1)}{4} (a+b)+x^2y^2(a\:b).
\end{align*}
\end{example}

\begin{proof}
It only remains to consider the surjectivity. If $x,y\neq 0$, for any composition $(a_1,\ldots,a_k)$, 
\[\kappa_{x,y}\left(
\begin{pmatrix}
a_1&0&\ldots&0\\
0&\ddots&\ddots&\vdots\\
\vdots&\ddots&\ddots&0\\
0&\ldots&0&a_k
\end{pmatrix}
\right)=x^ky^k (a_1,\ldots,a_k)+\mbox{a lin. comb. of compositions of length $<k$}.
\]
By a triangularity argument, as $x^ky^k \neq 0$, $\kappa_{x,y}$ is surjective.
If $x=0$, then for any $n\geq 1$, $H_n(x)=0$, so if $M\neq 1$, $\Psi_{0,y}(M)=0$, and $\Psi_{0,y}$ is the canonical projection on $\Q\subseteq \QSym$. 
The result is similar if $y=0$. 
\end{proof}

\begin{remark}
In particular, for any packed matrix $M$, as $H_k(1)=0$ for any $k\geq 2$,
\[\kappa_{1,1}(M)=\begin{cases}
\mathrm{diag}(M)\mbox{ if $M$ is diagonal},\\
0\mbox{ otherwise}.
\end{cases}\]
\end{remark}

\subsection{Two parameters shuffles}\label{ssec:sprod}

Recall that we have introduced in Definition \ref{bigrapack}
a bigraded bialgebra of packed matrices, $(\bfH_{\pack},\searrow,\bdeltares)$, with grading $\gr(M):=(\row(M), \col(M))$.
The dual bigraded bialgebra of packed matrices, $(\bfH_{\pack},\shuffle,\Delta)$ has a coproduct given by:

\begin{lemma}\label{bigraddual}
The product $\shuffle$ on $\bfH_\pack$ dual to $\bdeltares$ is given by
\begin{align*}
&\forall M,M'\in \pack,&M\shuffle M'&=\sum_{\substack{\sigma \in \sh(\row(M)\row(M'))\\ \tau \in \sh(\col(M),\col(M'))}} \mu(\sigma) (M\searrow M')\mu(\tau)^\top.
\end{align*}
\end{lemma}
\begin{proof}
The Lemma follows from Prop. \ref{propqshu}, keeping from the description of $\squplus$ only the components that preserve the graduation.
\end{proof}

\begin{example}
\begin{align*}
\begin{pmatrix}
1
\end{pmatrix}\shuffle \begin{pmatrix}
0&1\\1&0
\end{pmatrix}&=\begin{pmatrix}
1&0&0\\0&0&1\\0&1&0
\end{pmatrix}+2\begin{pmatrix}
0&1&0\\0&0&1\\1&0&0
\end{pmatrix}+3\begin{pmatrix}
0&0&1\\0&1&0\\1&0&0
\end{pmatrix}+2\begin{pmatrix}
0&0&1\\1&0&0\\0&1&0
\end{pmatrix}+\begin{pmatrix}
0&1&0\\1&0&0\\0&0&1
\end{pmatrix}.
\end{align*}
\end{example}

From the description of the Giusti-Lee-Nanda-Oberhauser two-parameter bialgebra structure on packed matrices \cite{Giusti} in \cite{Schmitz}, it follows that:

\begin{cor}
The bialgebra $(\bfH_{\pack},\shuffle,\Delta)$ identifies with the two-parameter bialgebra of Giusti-Lee-Nanda-Oberhauser, which is therefore by Lemma \ref{bigraddual} the bigraded dual bialgebra of  $(\bfH_{\pack},\searrow,\bdeltares)$.
\end{cor}

The following Proposition was obtained in \cite{Schmitz}. For completeness sake we include here a short, self-contained, proof and refer to the article for complementary results. Notation is as in the account of the classical Hoffman theorem provided in Section \ref{descents}.

\begin{prop}
The following is a coalgebra automorphism and bialgebra isomorphism:
\begin{align*}
\Upsilon&:\left\{\begin{array}{rcl}
(\bfH_\pack,\shuffle,\Delta)&\longrightarrow&(\bfH_\pack,\squplus,\Delta)\\
M\in \pack&\longmapsto&\displaystyle \sum_{\substack{\sigma \in \inc(\row(M))\\ \tau\in \inc(\col(M'))}} \frac{1}{\sigma!\tau!} \mu(\sigma) M\mu(\tau)^\top.
\end{array}\right.
\end{align*}
\end{prop}

\begin{proof}
For any packed matrix $M$, $\Upsilon(M)-M$ is in the linear span of matrices with strictly less rows or columns than $M$. A triangularity argument then shows that $\Upsilon$ is a bijection. 
By Corollary \ref{corbattages}, if $M$, $M'\in \pack$, with $k=\row(M)$, $k'=\row(M')$, $l=\col(M)$ and $l'=\col(M')$,
\begin{align*}
&\Upsilon(M\shuffle M')\\
&=\sum_{\substack{\sigma_1\in \sh(k,k')\\ \tau_1 \in \inc(k+k')}}\:\sum_{\substack{\sigma_2\in \sh(l,l')\\ \tau_2 \in \inc(l+l')}}
\frac{1}{\tau_1!\tau_2!} \mu(\tau_1\circ \sigma_1)(M\searrow M')\mu(\tau_2\circ \sigma_2)^\top\\
&=\sum_{\substack{\sigma'_1\in \inc(k)\\ \sigma''_1\in \inc(k')\\ \tau_1\in \qsh(\max(\sigma'_1),\max(\sigma''_1))}}\:
\sum_{\substack{\sigma'_2\in \inc(l)\\ \sigma''_2\in \inc(l')\\ \tau_2\in \qsh(\max(\sigma'_2),\max(\sigma''_2))}}
\frac{1}{\sigma'_1!\sigma''_1! \sigma'_2!\sigma''_2}\mu(\tau_1\circ (\sigma_1'\otimes \sigma_1''))(M\searrow M')\mu(\tau_2\circ (\sigma_2'\otimes \sigma_2''))^\top\\
&=\Upsilon(P)\squplus \Upsilon(M').
\end{align*} 
So $\Upsilon$ is an algebra morphism from $\shuffle$ to $\squplus$. \\

Let $M\in \pack$. We put $k=\row(M)$ and $l=\col(M)$.
\begin{align*}
\Delta \circ \Upsilon(M)&=\sum_{\substack{\sigma\in\inc(k)\\\tau\in \inc(l)}}\: \sum_{\mu(\sigma) M\mu(\tau)^\top=M'\searrow M''}\frac{1}{\sigma!\tau!}M'\otimes M''.
\end{align*}
If $\mu(\sigma) M\mu(\tau)=M'\searrow M''$, with $\sigma\in \inc(k)$ and $\tau\in \inc(l)$, then we can write un a unique way $M=M_1\searrow M_2$,
$\sigma=\sigma_1\otimes \sigma_2$ and $\tau=\tau_1\otimes \tau_2$, such that $M'=\mu(\sigma_1)M_1\mu(\tau_1)^\top$ and $M''=\mu(\sigma_2)M_1\mu(\tau_2)^\top$.
Therefore, 
\begin{align*}
\Delta \circ \Upsilon(M)&=\sum_{M=M_1\searrow M_2} \sum_{\substack{\sigma_1\in \inc(\row(M_1))\\
\sigma_2 \in \inc(\row(M_2))\\ \tau_1\in \inc(\col(M_1))\\ \tau_2\in \inc(\col(M_2))}}\frac{1}{\sigma_1!\sigma_2!\tau_1!\tau_2!}\mu(\sigma_1)M_1\mu(\tau_1)^\top
\otimes \mu(\sigma_2)M_2\mu(\tau_2)^\top\\
&=(\Upsilon \otimes \Upsilon)\circ \Delta(M).
\end{align*} 
Therefore, $\Upsilon$ is an automorphism of the coalgebra $(\bfH_\pack,\Delta)$. 
\end{proof}

\newcommand{\sym}{\mathfrak{S}}
\begin{remark}
The linear subspace generated by permutations matrices is a subbialgebra of $(\bfH_\pack,\shuffle,\Delta)$, but not of $(\bfH_\pack,\squplus,\Delta)$.
This induces a commutative bialgebraic structure on permutations, given by
\begin{align*}
&\forall \sigma \in \sym_k,\:\forall \tau\in \sym_l,&\sigma \shuffle \tau&=\sum_{\alpha,\beta \in \sh(k,l)} \alpha\circ (\sigma \otimes \tau)\circ \beta^{-1},\\
&\forall \sigma \in \sym_n,&\Delta(\sigma)&=\sum_{\sigma=\sigma'\otimes \sigma''}\sigma'\otimes \sigma''.
\end{align*}
This is not the product defined by Hivert, Novelli and Thibon in \cite{Hivert2008}, which we denote here by $\star$, although the coproduct is the same:
\begin{align*}
&\forall \sigma \in \sym_k,\:\forall \tau\in \sym_l,&\sigma \star \tau&=\sum_{\alpha\in \sh(k,l)} \alpha\circ (\sigma \otimes \tau)\circ \alpha^{-1}.
\end{align*}
For example,
\begin{align*}
(1)\shuffle(21)&=(132)+2(312)+3(321)+2(231)+(213),\\
(1)\star (21)&=(132)+(213)+(321).
\end{align*}
\end{remark}

\section{The double bialgebra structure}\label{sec:double}

\subsection{Polynomial realization}\label{ssec:const}
Let us recall first from \cite{Diehl} the construction of a polynomial realization of $\bfH_\pack$. 
\begin{defi}
A totally ordered alphabet is a pair $(X,\leq)$ such that $X$ is a finite or countable set and $\leq$ is a total order on $X$.
If $X$ and $Y$ are totally ordered alphabets, we denote by $\bfA_{X,Y}$ the algebra of formal power series 
\[\bfA_{X,Y}=\Q[[t_{i,j}\mid i\in X, j\in Y]].\]
\end{defi}

\begin{prop}\label{propinjectif}
Let $X,Y$ be two totally ordered alphabets.
For any $M=(m_{i,j})_{\substack{1\leq i\leq k\\1\leq j\leq l}}\in \pack$, we put
\[\Phi_M(X,Y)=\sum_{\substack{i_1<\ldots <i_k\:\mbox{\scriptsize in }X\\
{j_1<\ldots <j_l \:\mbox{\scriptsize in }Y}}} t_{i_1,j_1}^{m_{1,1}}\ldots t_{i_k,j_l}^{m_{k,l}}.\]
This is extended as a linear map $\Phi_\bullet(X,Y):\bfH_\pack\longrightarrow \bfA_{X,Y}$. It is injective if, and only if, $X$ and $Y$ are infinite.
\end{prop}

\begin{example}
Let $a,b\in\N_{\geq 1}$. 
\begin{align*}
\Phi_{\begin{pmatrix}
a
\end{pmatrix}}(X,Y)&=\sum_{i\in X,\: j\in Y}t_{i,j}^a,\\
\Phi_{\begin{pmatrix}
a&b
\end{pmatrix}}(X,Y)&=\sum_{i\in X,\:j_1<j_2\in Y}t_{i,j_1}^at_{i,j_2}^b,\\
\Phi_{\begin{pmatrix}
a\\b
\end{pmatrix}}(X,Y)&=\sum_{i_1<i_2\in X,\:j\in Y}t_{i_1,j}^at_{i_2,j}^b,\\
\Phi_{\begin{pmatrix}
a&0\\0&b
\end{pmatrix}}(X,Y)&=\sum_{i_1<i_2\in X,\:j_1<j_2\in Y}t_{i_1,j_1}^at_{i_2,j_2}^b,\\
\Phi_{\begin{pmatrix}
0&b\\a&0
\end{pmatrix}}(X,Y)&=\sum_{i_1<i_2\in X,j_1<j_2\in Y}t_{i_1,j_2}^bt_{i_2,j_1}^a.
\end{align*}
\end{example}

\begin{proof}
Let us denote by $\bfM_{X,Y}$ the set of monomials of $\bfA_{X,Y}$.
If $P=t_{i_1,j_1}^{a_1}\ldots t_{i_k,j_k}^{a_k}\in \bfM_{X,Y}$, with $a_1,\ldots,a_k\geq 1$, we put:
\begin{itemize}
\item $\{i_1,\ldots,i_k\}=\{i'_1<\ldots<i'_p\}$, totally ordered subset of $X$.
\item $\{j_1,\ldots,j_k\}=\{j'_1<\ldots<j'_q\}$, totally ordered subset of $Y$.
\item For any $r\in [p]$, for any $s\in [q]$, the power of $t_{i'_r,j'_s}$ in $P$ is denoted by $m_{r,s}$.
\end{itemize}
We then put
\[\varsigma(P)=(m_{r,s})_{\substack{1\leq r\leq p\\1\leq s\leq q}}.\]
For any $r\in [p]$, by construction of $\{i'_1<\ldots<i'_p\}$, there exists $r'$ such $i_{r'}=i'_r$. Let $s\in [q]$ such that $j_{r'}=j'_s$.
Then $m_{r,s}=a_{r'}\neq 0$: $\varsigma(P)$ has no zero row. Similarly, it has no zero column, so it belongs to $\pack$. Its weight is the degree of $P$. 
Then, by construction of $\Phi_M(X,Y)$,
\[\Phi_M(X,Y)=\sum_{P\in \bfM_{X,Y},\: \varsigma(P)=M} P.\]
This implies that the kernel of $\Phi_M(X,Y)$ is generated by the packed matrices $M$ such that $\Phi_M(X,Y)=0$, that is to say such that there is no $P\in \bfM_{X,Y}$
such that $\varsigma(P)=M$. If $X$ or $Y$ is finite, then taking a packed matrix $M$ such that the number of rows is strictly greater than the cardinality of $X$,
or such that  the number of columns is strictly greater than the cardinality of $Y$, then $\Phi_M(X,Y)=0$, and $\Phi_\bullet(X,Y)$ is not injective.
If $X$ and $Y$ are infinite, for any $M\in \pack$, there exists a monomial $P\in \bfM_{X,Y}$ such that $\varsigma(P)=M$,
so $\Phi_\bullet(X,Y)$ is injective. 
\end{proof}

\begin{prop}
Let $M,M'\in \pack$, with respectively $k$ and $k'$ rows, and $l$ and $l'$ columns. 
Then, for any totally ordered alphabets $X,Y$,
\[\Phi_M(X,Y)\Phi_{M'}(X,Y)=\Phi_{M\squplus M'}(X,Y).\]
\end{prop}

\begin{proof}
Let $P,P'\in \bfM_{X,Y}$ such that $\varsigma(P)=M$ and $\varsigma(P')=M'$. Then $MM'$ appears in $\Phi_M(X,Y)\Phi_{M'}(X,Y)$.
Let $i_1<\ldots<i_k$ in $X$ the left indices of the indeterminates appearing in $M$ and $i'_1<\ldots<i'_{k'}$ the left indices of the indeterminates appearing in $M'$.
We then put $\{i_1,\ldots,i_k,i'_1,\ldots,i'_{k'}\}=\{i''_1<\ldots<i''_{k''}\}$, totally ordered subset of $X$. This defines a map $\sigma:[k+k']\longrightarrow [k'']$, such that:
\begin{itemize}
\item For any $p\leq k$, $i''_{\sigma(p)}=i_p$.
\item For any $p>k$, $i''_{\sigma(p)}=i'_{p-k}$. 
\end{itemize} 
As $i_1<\ldots<i_k$ in $X$, necessarily $i''_{\sigma(1)}<\ldots i''_{\sigma(k)}$ in $X,$ so $\sigma(1)<\ldots<\sigma(k)$.
Similarly, $\sigma(k+1)<\ldots<\sigma(k+k')$, and $\sigma \in \qsh(k,k')$.
Considering now the right indices of the indeterminates appearing in $M$ and $M'$, we obtain a $(l,l')$-quasi-shuffle $\tau$, such that
\[\varsigma(PP')=\mu(\sigma)(M\searrow M')\mu(\tau)^\top.\]
Conversely, if $P$ is a monomial appearing in $\Phi_{\mu(\sigma)(M\searrow M')\mu(\tau)^\top}(X,Y)$ for given $(k,k')$- and $(l,l')$-quasi-shuffles,
then $P$ is the product of two monomials $Q$ and $Q'$, appearing respectively in $\Phi_M(X,Y)$ and $\Phi_{M'}(X,Y)$.
So $\Phi_M(X,Y)\Phi_{M'}(X,Y)=\Phi_{M\squplus M'}(X,Y)$. 
\end{proof}

\begin{notation}
Let $X,Y$ be two totally ordered alphabets. The set $X\sqcup Y$,
disjoint union of $X$ and $Y$, is given a total order by preserving the total orders on $X$ and $Y$ and assigning that $x<y$ for any $(x,y)\in X\times Y$.
This total ordered alphabet is denoted by $X+Y$. Note that $X+Y$ and $Y+X$ are different in general, even if the underlying sets are the sames.
If $X_1,X_2,Y_1,Y_2$ are four totally ordered alphabets, we consider the algebra map
\[\pi_{(X_1,X_2),(Y_1,Y_2)}:\left\{\begin{array}{rcl}
\bfA_{X_1+X_2,Y_1+Y_2}&\longrightarrow&\bfA_{X_1,Y_1}\otimes \bfA_{X_2,Y_2}\\
t_{i,j}&\longmapsto&\begin{cases}
t_{i,j}\otimes 1\mbox{ if }(i,j)\in X_1\times Y_1,\\
1\otimes t_{i,j}\mbox{ if }(i,j)\in X_2\times Y_2,\\
0\mbox{ otherwise}.
\end{cases}\end{array}\right.\]
\end{notation}

\begin{prop}
Let $X_1,X_2,Y_1,Y_2$ be four totally ordered alphabets. For any $M\in \pack$,
\begin{align*}
\pi_{(X_1,X_2),(Y_1,Y_2)}(\Phi_M(X_1+X_2,Y_1+Y_2))&=\sum_{M=M_1\searrow M_2}\Phi_{M_1}(X_1,Y_1)\otimes \Phi_{M_2}(X_2,Y_2)\\
&=(\Phi_\bullet(X)\otimes \Phi_\bullet(Y))\circ \Delta(M).
\end{align*}
\end{prop}

\begin{proof}
Let $P\in \calM_{X_1+X_2,Y_1+Y_2}$ such that $\varsigma(P)=M$. If $\pi_{(X_1,X_2),(Y_1,Y_2)}(P)\neq 0$, then $P$ can be written as a product $P_1P_2$,
with $P_1\in \calM_{X_1,Y_1}$ and $P_2\in \calM_{X_2,Y_2}$. By definition of the sum of two ordered alphabets, $M=\varsigma(P)=\varsigma(P_1)\searrow \varsigma(P_2)$. 
Conversely, if $M=M_1\searrow M_2$, for any $P_1\in \calM_{X_1,Y_1}$ such that $\varsigma(P_1)=M_1$ and  $P_2\in \calM_{X_2,Y_2}$ such that $\varsigma(P_2)=M_2$,
then $\varsigma(P_1P_2)=\varsigma(P_1)\searrow \varsigma(P_2)=M$ and $\pi_{(X_1,X_2),(Y_1,Y_2)}(P_1P_2)=P_1\otimes P_2$. This implies the announced formula.
\end{proof}

\subsection{The second coproduct and the double bialgebra structure}

We face now the key problem, that we had postponed, of defining a double bialgebra structure on $(\bfH_\pack,\squplus,\Delta)$ by introducing another coproduct $\delta$ with suitable compatibility properties. By duality this new coproduct would also define still another algebra (and bialgebra) structure on $\bfH_\pack$ with compatibility properties with respect to $\searrow$ and $\blacktriangle$, as expected from the case of classical noncommutative symmetric functions.
We use the polynomial realization of $\bfH_\pack$ and the trick of considering products of alphabets to achieve the construction.

\begin{notation}
Let $X$ and $Y$ be two totally ordered alphabets. Then $X\times Y$ is totally ordered by the lexicographic order:
\begin{align*}
&\forall (x,y),(x',y')\in X\times Y,&(x,y)\leq (x',y')&\mbox{ if }(x<x')\mbox{ or }(x=x'\mbox{ and }y\leq y').
\end{align*}
This totally ordered alphabet is denoted by $XY$. For any totally ordered alphabets $X_1$, $X_2,$ $Y_1$ and $Y_2$, we consider the algebra map
\[\iota_{(X_1,X_2),(Y_1,Y_2)}:\left\{\begin{array}{rcl}
\bfA_{X_1X_2,Y_1Y_2}&\longrightarrow&\bfA_{X_1,Y_1}\otimes \bfA_{X_2,Y_2}\\
t_{(x_1,x_2),(y_1,y_2)}&\longmapsto&t_{x_1,y_1}\otimes t_{x_2,y_2}.
\end{array}\right.\]
\end{notation}

\begin{prop}
Let $k\in \N$. We denote by $\adm(k)$ the set of pairs $(\sigma',\sigma'')$ such that:
\begin{itemize}
\item $\sigma':[k]\longrightarrow [k']$ and $\sigma'':[k]\longrightarrow [k'']$ are surjective maps.
\item$\sigma'\in \inc(k)$.
\item For any $i,j\in [k]$, ($i<j$ and $\sigma'(i)=\sigma'(j)$) $\Longrightarrow$ $\sigma''(i)<\sigma''(j)$.
\end{itemize}
Let $X_1$, $X_2$, $Y_1$ and $Y_2$ be totally ordered alphabets. For any $M\in \pack$, 
\begin{align*}
\iota_{(X_1,X_2),(Y_1,Y_2)}\circ \Phi_M(X_1X_2,Y_1Y_2)&=\sum_{\substack{(\sigma',\sigma'')\in \adm(\row(M))\\ (\tau',\tau'')\in \adm(\col(M))}}
\Phi_{\mu(\sigma')M\mu(\tau')^\top}(X_1,Y_1)\otimes \Phi_{\mu(\sigma'')M\mu(\tau'')^\top}(X_2,Y_2).
\end{align*}\end{prop}

\begin{proof}
Let $P\in \calM_{X_1X_2,Y_1Y_2}$ such that $\varsigma(P)=M$. Let us denote by $(i_1,j_1)<\ldots<(i_k,j_k)\in X_1X_2$ the left indices of the indeterminates
which appears in $P$. Let us put 
\[\{i_1,\ldots,i_k\}=\{i'_1<\ldots<i'_{k'}\}\subseteq X_1.\]
This defines a surjective map $\sigma':[k]\longrightarrow [k']$ sending any $p\in [k]$ to the unique $\sigma'(p)\in[k']$ such that $i'_{\sigma'(p)}=i_p$. 
By definition of the lexicographic order, $i_1\leq \ldots \leq i_k$, so $\sigma'(1)\leq \ldots \leq\sigma'(k)$. 
Let us put 
\[\{j_1,\ldots,j_k\}=\{j''_1<\ldots<j''_{k''}\}\subseteq X_2.\]
This defines a surjective map $\sigma'':[k]\longrightarrow [k'']$
sending any $q\in [k]$ on the unique $\sigma''(q)\in[k'']$ such that $j''_{\sigma''(q)}=j_q$. If $p<q$ and $\sigma'(p)=\sigma'(q)$, then $i_p=i_q$.
By definition of the lexicographic order, $j_p<j_q$, so $\sigma''(p)<\sigma''(q)$. We proved that $(\sigma',\sigma'')\in \adm(k)$.
Working with the right indices, we similarly obtain $(\tau,\tau')\in \adm(l)$, such that 
\[(\varsigma\otimes \varsigma)\circ \iota_{(X_1,X_2),(Y_1,Y_2)}(P)=\mu(\sigma')M\mu(\tau')^\top \otimes \mu(\sigma'')M\mu(\tau'')^\top.\]
Conversely, let $(\sigma',\sigma'')\in \adm(k)$, $(\tau',\tau'')\in \adm(l)$, $P_1\in \bfA_{X_1,Y_1}$ such that $\varsigma(P_1)=\mu(\sigma')M\mu(\tau')^\top$
and $P_2\in \bfA_{X_2,Y_2}$ such that $\varsigma(P_2)=\mu(\sigma'')M\mu(\tau'')^\top$.
 Each indeterminate in $P_1$ and $P_2$ corresponds to an entry of $M$; joining together the left and right indices of the indeterminates in $P_1$ and $P_2$
corresponding to the same entry of $M$, we obtain $P\in \bfA_{X_1X_2,Y_1Y_2}$ such that $\iota_{(X_1,X_2),(Y_1,Y_2)}(P)=P_1\otimes P_2$. The definition of $\adm(k)$ and $\adm(l)$
yields $\varsigma(P)=M$. This concludes the proof.    
\end{proof}

\begin{theo}
We define a second coproduct $\delta$ on $\bfH_\pack$, sending any packed matrix $M$ with $k$ rows and $l$ columns on
\begin{align*}
\delta(M)&=\sum_{\substack{(\sigma',\sigma'')\in \adm(k)\\ (\tau',\tau'')\in \adm(l)}}
\mu(\sigma')M\mu(\tau')^\top\otimes\mu(\sigma'')M\mu(\tau'')^\top.
\end{align*}
Then $(\bfH_\pack,\squplus,\Delta,\delta)$ is a double bialgebra.
\end{theo}

\begin{proof}
Let us choose six infinite totally ordered alphabets $X_1,X_2,X_3,Y_1,Y_2,Y_3$.
Observe that for any totally ordered alphabets $X$, $Y$ and $Z$, 
\[(XY)Z=X(YZ).\] Hence,
\begin{align*}
&(\Phi_\bullet(X_1,Y_1)\otimes \Phi_\bullet(X_2,Y_2)\otimes \Phi_\bullet(X_3,Y_3))\circ (\delta \otimes \id)\circ \delta\\ 
&=(\iota_{(X_1,X_2),(Y_1,Y_2)}\otimes \id)\circ \iota_{(X_1X_2,X_3),(Y_1Y_2,Y_3)}\circ \Phi_\bullet(X_1X_2X_3,Y_1Y_2Y_3)\\
&=(\id \otimes \iota_{(X_2,X_3),(Y_2,Y_3)})\circ \iota_{(X_1,X_2X_3),(Y_1,Y_2Y_3)}\circ \Phi_\bullet(X_1X_2X_3,Y_1Y_2Y_3)\\
&=(\Phi_\bullet(X_1,Y_1)\otimes \Phi_\bullet(X_2,Y_2)\otimes \Phi_\bullet(X_3,Y_3))\circ (\id \otimes \delta)\circ \delta. 
\end{align*} 
As $\Phi_\bullet(X_i,Y_i)$ is injective for $i\in [3]$ (Proposition \ref{propinjectif}), we conclude that $\delta$ is coassociative.\\

Let $M,M'\in \pack$.
\begin{align*}
&(\Phi_\bullet(X_1,Y_1)\otimes \Phi_\bullet(X_2,Y_2))\circ \delta(M\squplus M')\\
&=\iota_{(X_1,X_2),(Y_1,Y_2)}(\Phi_{M\squplus M'}(X_1X_2,Y_1Y_2))\\
&=\iota_{(X_1,X_2),(Y_1,Y_2)}(\Phi_M(X_1X_2,Y_1Y_2)\Phi_{M'}(X_1X_2,Y_1Y_2))\\
&=\iota_{(X_1,X_2),(Y_1,Y_2)}(\Phi_M(X_1X_2,Y_1Y_2))\iota_{(X_1,X_2),(Y_1,Y_2)}(\Phi_{M'}(X_1X_2,Y_1Y_2))\\
&=(\Phi_\bullet(X_1,Y_1)\otimes \Phi_\bullet(X_2,Y_2))(\delta(M)\squplus\delta(M')).
\end{align*}
Again,the injectivity of $\Phi_\bullet(X_i,Y_i)$ for $i\in [2]$ implies that $\delta$ is multiplicative.   \\

Note that for any totally ordered alphabets, $(X+Y)Z=XZ+YZ$
(but $X(Y+Z)$ and $XY+XZ$ may be different). Let $M$ be a packed matrix and write $m_{1,3,24}(a\otimes b\otimes c\otimes d):=a\otimes c\otimes bd$, then:
\begin{align*}
&(\Phi_\bullet(X_1,X_2)\otimes \Phi_\bullet(Y_1,Y_2)\otimes \Phi_\bullet(X_3,Y_3))\circ (\Delta \otimes \id)\circ \delta(M)\\
&=(\pi_{(X_1,X_2),(Y_1,Y_2)} \otimes \id)\circ \iota_{(X_1+X_2,X_3),(Y_1+Y_2,Y_3)}(\Phi_M((X_1+X_2)X_3,(Y_1+Y_2)Y_3) \\
&=m_{1,3,24}(\iota_{(X_1,X_3),(Y_1,Y_3)} \otimes \iota_{(X_2,X_3),(Y_2,Y_3)})\circ \pi_{(X_1X_3,X_2X_3),(Y_1Y_3,Y_2Y_3)}
(\Phi_M(X_1X_3+X_2X_3,Y_1Y_3+Y_2Y_3))\\
&=(\Phi_\bullet(X_1,X_2)\otimes \Phi_\bullet(Y_1,Y_2)\otimes \Phi_\bullet(X_3,Y_3))\circ \squplus_{1,3,24}\circ (\delta \otimes \delta)\circ \Delta(M).
\end{align*}
Again, an injectivity argument proves the compatibility between the two coproducts.\\

Let us conclude on the counit of $\delta$. 
Let us define $\epsilon_\delta$ by
\begin{align*}
&\forall M\in \pack,&\epsilon_\delta(M)&=\begin{cases}
1\mbox{ if }M=1,\\
1\mbox{ if $M$ has one row and one column},\\
0\mbox{ otherwise}.
\end{cases}
\end{align*}
Let $M$ be a nonempty packed matrix. Let $(\sigma',\sigma'')\in \adm(\row(M))$, such that $\sigma'([\row(M)])=[1]$.
By definition of $\adm(\row(M))$, $\sigma''=\id$. Therefore,
\begin{align*}
(\epsilon_\delta \otimes \id)\circ \delta(M)&=\epsilon_\delta\left(\begin{pmatrix}
1&\ldots&1
\end{pmatrix}
M \begin{pmatrix}
1\\ \vdots\\1
\end{pmatrix}\right)M+0=M.
\end{align*}
Similarly, if $\sigma''([\row(M)])=[1]$, necessarily $\sigma'$ is injective. As $\sigma'\in \inc(k)$, it is the identity, so 
\begin{align*}
(\id \otimes\epsilon_\delta)\circ \delta(M)&=\epsilon_\delta\left(\begin{pmatrix}
1&\ldots&1
\end{pmatrix}
M \begin{pmatrix}
1\\ \vdots\\1
\end{pmatrix}\right)M+0=M.
\end{align*}
So $\epsilon_\delta$ is the counit of $\delta$.
\end{proof}

\begin{example} 
\begin{align*}
\delta\left(\begin{pmatrix}
a
\end{pmatrix}\right)&=\begin{pmatrix}
a
\end{pmatrix}\otimes \begin{pmatrix}
a
\end{pmatrix},\\
\delta\left(\begin{pmatrix}
a&b
\end{pmatrix}\right)&=\begin{pmatrix}
a+b
\end{pmatrix}\otimes \begin{pmatrix}
a&b
\end{pmatrix}+\begin{pmatrix}
a&b
\end{pmatrix}\otimes\left(\begin{pmatrix}
a&b
\end{pmatrix}+\begin{pmatrix}
b&a
\end{pmatrix}+\begin{pmatrix}
a+b
\end{pmatrix}\right),\\
\delta\left(\begin{pmatrix}
a\\b
\end{pmatrix}\right)&=\begin{pmatrix}
a+b
\end{pmatrix}\otimes \begin{pmatrix}
a\\b
\end{pmatrix}+\begin{pmatrix}
a\\b
\end{pmatrix}\otimes\left(\begin{pmatrix}
a\\b
\end{pmatrix}+\begin{pmatrix}
b\\a
\end{pmatrix}+\begin{pmatrix}
a+b
\end{pmatrix}\right),\\
\delta\left(\begin{pmatrix}
a&b\\c&d
\end{pmatrix}\right)&=\begin{pmatrix}
a+b+c+d
\end{pmatrix}\otimes \begin{pmatrix}
a&b\\c&d
\end{pmatrix}+\begin{pmatrix}
a+c&b+d
\end{pmatrix}\otimes \left(\begin{pmatrix}
a&b\\c&d
\end{pmatrix}+\begin{pmatrix}
b&a\\d&c
\end{pmatrix}+\begin{pmatrix}
a+b\\c+d
\end{pmatrix}\right)\\
&+\begin{pmatrix}
a+b\\c+d
\end{pmatrix}\otimes\left(\begin{pmatrix}
a&b\\c&d
\end{pmatrix}+\begin{pmatrix}
c&d\\a&b
\end{pmatrix}+\begin{pmatrix}
a+c&b+d
\end{pmatrix}\right)\\
&+\begin{pmatrix}
a&b\\c&d
\end{pmatrix}\otimes
 \left(\begin{array}{c}
\begin{pmatrix}
a&b\\c&d
\end{pmatrix}+\begin{pmatrix}
c&d\\a&b
\end{pmatrix}+\begin{pmatrix}
a+c&b+d
\end{pmatrix}+\begin{pmatrix}
b&a\\d&c
\end{pmatrix}+\begin{pmatrix}
d&c\\b&a
\end{pmatrix}\\
+\begin{pmatrix}
b+d&a+c
\end{pmatrix}+\begin{pmatrix}
a+b\\c+d
\end{pmatrix}+\begin{pmatrix}
c+d\\a+b
\end{pmatrix}+\begin{pmatrix}
a+b+c+d
\end{pmatrix}
\end{array}\right).
\end{align*}
\end{example}

\subsection{Double bialgebra morphisms}

\begin{notation}
Let $X,Y$ be two totally ordered alphabets. We denote by $\tau$ the algebra morphism defined by
\[\tau_{X,Y}:\left\{\begin{array}{rcl}
\bfA_{X,Y}&\longrightarrow&\bfA_{Y,X}\\
t_{x,y}&\longmapsto&t_{y,x}.
\end{array}\right.\]
\end{notation}

Recall the map $T$  of Proposition \ref{propT2}:
\[T:\left\{\begin{array}{rcl}
\bfH_\pack&\longrightarrow&\bfH_\pack\\
M\in \pack&\longmapsto&M^\top.
\end{array}
\right.\]

\begin{prop}
For any totally ordered alphabets $X,Y$,
\[\tau_{X,Y}\circ \Phi_\bullet(X,Y)=\Phi_\bullet(Y,X)\circ T. \]
Therefore, the map $T$  is a double bialgebra automorphism. 
\end{prop}

\begin{proof}
We already know that $T$ is a bialgebra automorphism of $(\bfH_\pack,\squplus,\Delta)$
Let $X,Y$ be two totally ordered alphabets. For any  $M=(m_{i,j})_{\substack{1\leq i\leq k\\ 1\leq j\leq l}}\in \pack$,
\begin{align*}
\tau_{X,Y}\circ \Phi_M(X,Y)&=\sum_{\substack{i_1<\ldots<i_k\: \mbox{\scriptsize in $X$},\\{j_1<\ldots<j_l\: \mbox{\scriptsize in $Y$}}}} t_{j_1,i_1}^{m_{1,1}}\ldots t_{j_l,i_k}^{m_{i_k,j_l}}=\Phi_{M^\top}(Y,X),
\end{align*}
so indeed $\tau_{X,Y}\circ \Phi_\bullet(X,Y)=\Phi_\bullet(Y,X)\circ T$. \\

Let $X_1,X_2,Y_1,Y_2$ be two infinite totally ordered alphabets.
\begin{align*}
(\Phi_\bullet(Y_1,X_1)\otimes \Phi_\bullet(Y_2,X_2))\circ (T\otimes T)\circ \delta
&=(\tau_{X_1,Y_1} \otimes \tau_{X_2,Y_2})\circ (\Phi_\bullet(X_1,Y_1)\otimes \Phi_\bullet(X_2,Y_2))\circ \delta\\
&=(\tau_{X_1,Y_1} \otimes \tau_{X_2,Y_2})\circ \iota_{(X_1,X_2),(Y_1,Y_2)}\circ \Phi_\bullet(X_1X_2,Y_1Y_2)\\
&=\iota_{(Y_1,Y_2),(X_1,X_2)}\circ \tau_{X_1X_2,Y_1Y_2} \circ \Phi_\bullet(X_1X_2,Y_1Y_2)\\
&=\iota_{(Y_1,Y_2),(X_1,X_2)}\circ \Phi_\bullet(Y_1Y_2,X_1X_2)\circ T\\
&=(\Phi_\bullet(Y_1,X_1)\otimes \Phi_\bullet(Y_2,X_2))\circ \delta \circ T.
\end{align*}
The injectivity of $\Phi_\bullet(Y_1,X_1)\otimes \Phi_\bullet(Y_2,X_2)$ implies that $(T\otimes T)\circ \delta=\delta \circ T$.
So $T$ is a double bialgebra automorphism. 
\end{proof}

\begin{notation}
\begin{enumerate}
\item If $X$ is a totally ordered alphabet, we denote by $\bfA_X$ the algebra of formal series
\[\bfA_X=\Q[[t_i\mid i\in X]].\]

Note that $\bfA_{X,Y}=\bfA_{XY}$.
\item Let us now return on the polynomial realization of $\QSym$ briefly discussed already in Section \ref{descents}.
We associate to  any composition $(a_1,\ldots,a_k)$ a (generalized) polynomial (corresponding to the usual  monomial quasi-symmetric function in the sense if Gessel \cite{Gessel1984}) by
\[\phi_{(a_1,\ldots,a_k)}(X)=\sum_{i_1<\ldots<i_k\:\mbox{\scriptsize in }X}t_{i_1}^{a_1}\ldots t_{i_k}^{a_k}.\]
This defines an algebra morphism $\phi_\bullet(X):\QSym\longrightarrow \bfA_X$, injective if and only if $X$ is infinite.
Moreover, if $X$ and $Y$ are two totally ordered alphabets,
\begin{align*}
(\phi_\bullet(X)\otimes \phi_\bullet(Y))\circ \Delta&=\pi_{X,Y}\circ \phi_\bullet(X+Y),\\
(\phi_\bullet(X)\otimes \phi_\bullet(Y))\circ \delta&=\iota_{X,Y}\circ \phi_\bullet(XY),
\end{align*} 
where $\pi_{X,Y}$ and $\iota_{X,Y}$ are the continuous algebra maps defined by
\begin{align*}
\pi_{X,Y}&:\left\{\begin{array}{rcl}
\bfA_{X+Y}&\longrightarrow&\bfA_X\otimes \bfA_Y\\
t_i&\longmapsto&\begin{cases}
t_i\otimes 1\mbox{ if }i\in X,\\
1\otimes t_i\mbox{ if }i\in Y,
\end{cases}
\end{array}\right.&
\iota_{X,Y}&:\left\{\begin{array}{rcl}
\bfA_{XY}&\longrightarrow&\bfA_X\otimes \bfA_Y\\
t_{(i,j)}&\longmapsto&t_i\otimes t_j.
\end{array}\right.
\end{align*}
\end{enumerate}
\end{notation}

\begin{prop}
The map $\theta:\QSym\longrightarrow \bfH_\pack$ of Proposition \ref{proptheta2} is a double bialgebra morphism.
\end{prop}

\begin{proof}
We already know that $\theta$ is a bialgebra morphism from $(\QSym,\squplus,\Delta)$ to $(\bfH_\pack,\squplus,\Delta)$. 
Let $X,Y$ be two totally ordered alphabets. Let us first show that
\[\phi_\bullet(XY)=\Phi_\bullet(X,Y)\circ \theta.\]
Indeed, if $(a_1,\ldots,a_k)$ is a composition and $P$ is a monomial appearing in $\phi_{(a_1,\ldots,a_k)}(XY)$, then there exists a unique packed matrix $M$
such that $P$ appears in $\Phi_M(X,Y)$; moreover, $\comp(M)=(a_1,\ldots,a_k)$. Conversely, if $\comp(M)=(a_1,\ldots,a_k)$ 
and $P$ is a monomial appearing in $\phi_M(X,Y)$, then $P$ also appears in $\phi_{(a_1,\ldots,a_k)}(XY)$, by definition of the lexicographic order.\\

Thus
\begin{align*}
(\Phi_\bullet(X_1,Y_1)\otimes \Phi_\bullet(X_2,Y_2))\circ (\theta \otimes \theta)\circ \delta
&=(\phi_\bullet(X_1X_2)\otimes \phi_\bullet(Y_1Y_2))\circ \delta\\
&=\iota_{X_1X_2,Y_1Y_2}\circ \phi_\bullet((X_1X_2)(Y_1Y_2))\\
&=\iota_{X_1X_2,Y_1Y_2}\circ \Phi_\bullet(X_1X_2,Y_1Y_2)\circ \theta\\
&=(\Phi_\bullet(X_1,Y_1)\otimes \Phi_\bullet(X_2,Y_2))\circ \delta\circ \theta.
\end{align*}
The injectivity of $\Phi_\bullet(X_i,Y_i)$ for $i\in [2]$ implies the compatibility of $\theta$ with $\delta$.
\end{proof}

\begin{theo}
Let $x,y\in \Q$. Then the morphism  $\kappa_{x,y}$  from $\bfH_\pack$ to $\QSym$ of Proposition \ref{propkxy2} is a double bialgebra morphism if, and only if, $x=y=1$.
\end{theo}

\begin{proof}
$\Longrightarrow$. Under the hypothesis, one has $\epsilon_\delta \circ \kappa_{x,y}=\epsilon_\delta$. Therefore,
\[1=\epsilon_\delta\circ \kappa_{x,y}\begin{pmatrix}
1
\end{pmatrix}=xy\epsilon_\delta((1))=xy,\]
so $xy=1$. Moreover,
\[0=\epsilon_\delta\circ \kappa_{x,y}\begin{pmatrix}
1&1
\end{pmatrix}=\frac{xy(y-1)}{2}\epsilon_\delta((2))=\frac{y-1}{2},\]
so $y=1$ and finally $x=y=1$. \\

$\Longleftarrow$. The map $\kappa_{1,1}$ is obviously compatible with $\epsilon_\delta$. For any totally ordered alphabet $X$, let us consider the algebra map
\[d_X:\left\{\begin{array}{rcl}
\bfA_{X,Y}&\longrightarrow&\bfA_X\\
t_{i,j}&\longmapsto&\begin{cases}
t_i\mbox{ if }i=j,\\
0\mbox{ otherwise}.
\end{cases}
\end{array}\right.\]
Let us first prove that $\phi_\bullet(X)\circ \kappa_{1,1}=d_X\circ\Phi_\bullet(X,X)$.
Let $M$ be a packed matrix. Then
\[\phi_\bullet(X)\circ \kappa_{1,1}(M)=\begin{cases}
\displaystyle \sum_{i_1<\ldots<i_{\row(M)} \:\mbox{\scriptsize in } X} t_{i_1}^{m_{1,1}}\ldots t_{i_{\row(M)}}^{m_{\row(M),\row(M)}}\mbox{ if $M$ is diagonal},\\
0\mbox{ otherwise}. 
\end{cases}\]

If $M$ is diagonal, we immediately obtain that the only monomials appearing in $\Phi_M(X,X)$ which do not vanish under $d_X$ are the monomials
\[t_{i_1,i_1}^{m_{1,1}}\ldots t_{i_{\row(M)},i_{\row(M)}}^{m_{\row(M),\row(M)}},\]
with $i_1<\ldots<i_{\row(M)}$ in $X$. In this case, we indeed have 
$\phi_\bullet(X)\circ \kappa_{1,1}(M)=d_X(\Phi_\bullet(X,X))$.
 
If $d_X(\Phi_M(X,X))\neq 0$, then in $\Phi_M(X,X)$ appears a monomial of the form $t_{i_1,i_1}^{m_1}\ldots t_{i_k,i_k}^{m_k}$, with $i_1<\ldots<i_k$ in $X$ and $m_1,\ldots,m_k\geq 1$. Hence, $M=\mathrm{diag}[m_1,\ldots,m_k]$ is diagonal. By contraposition, if $M$ is not diagonal, then $d_X(\Phi_M(X,X))=0=\phi_\bullet(X)\circ \kappa_{1,1}(M)$. \\

Let us now choose two infinite totally ordered alphabets $X$ and $Y$. Then
\begin{align*}
(\phi_\bullet(X)\otimes \phi_\bullet(Y))\circ \delta \circ \kappa_{1,1}&=\iota_{X,Y}\circ \phi_\bullet(XY)\circ \kappa_{1,1}\\
&=\iota_{X,Y}\circ d_{XY}\circ \Phi_\bullet(XY,XY)\\
&=(d_X\otimes d_Y)\circ \iota_{(X,Y),(X,Y)}\circ \Phi_\bullet(XY,XY)\\
&=(d_X\otimes d_Y)\circ (\Phi_\bullet(X,X)\otimes \Phi_\bullet(Y,Y))\circ \delta\\
&=(\phi_\bullet(X)\otimes \phi_\bullet(Y))\circ (\kappa_{1,1}\otimes \kappa_{1,1})\circ \delta.
\end{align*}
We conclude with the injectivity of $\phi_\bullet(X)$ and $\phi_\bullet(Y)$ that $ \delta \circ \kappa_{1,1}=(\kappa_{1,1}\otimes \kappa_{1,1})\circ \delta$.
So $\kappa_{1,1}$ is a double bialgebra morphism. 
\end{proof}

\begin{cor}
The unique double bialgebra morphism $\phi_{\bfH_\pack}$ from $(\bfH_\pack,\squplus,\Delta,\delta)$ to $(\K[X],m,\Delta,\delta)$ is given by
\begin{align*}
&\forall M\in \pack,&\phi_{\bfH_\pack}(M)&=\begin{cases}
H_{\row(M)}(X)\mbox{ if $M$ is diagonal},\\
0\mbox{ otherwise}.
\end{cases}
\end{align*}
\end{cor}

\begin{proof}
By composition, $\phi_{\QSym}\circ \kappa_{1,1}$ is a double bialgebra morphism, so by unicity is equal to $\phi_{\bfH_\pack}$. For any packed matrix $M$,
\begin{align*}
\phi_{\bfH_\pack}(M)=\phi_{\QSym}\circ \kappa_{1,1}(M)&=\begin{cases}
\phi_{\QSym}(\mathrm{diag}(M))\mbox{ if $M$ is diagonal},\\
0\mbox{ otherwise}
\end{cases}\\
&=\begin{cases}
H_{\row(M)}(X)\mbox{ if $M$ is diagonal},\\
0\mbox{ otherwise}.
\end{cases} \qedhere
\end{align*}
\end{proof}

\small
\vspace{1cm}

\noindent {\bf Acknowledgements.} 
The authors acknowledge support from the grant ANR-20-CE40-0007
Combinatoire Alg\'ebrique, R\'esurgence, Probabilit\'es Libres et Op\'erades.
L.F and F.P. would also like to thank Sapienza Universit\`a di Roma for its hospitality.
F.P. was also supported by the ANR -- FWF project PAGCAP.

\bibliographystyle{amsplain}
\bibliography{biblio}

\noindent
Lo\"\i c Foissy\\
Laboratoire de Math\'ematiques Pures et Appliqu\'ees Joseph Liouville\\
Universit\'e du Littoral C\^ote d'opale\\ Centre Universitaire de la Mi-Voix\\ 50, rue Ferdinand Buisson, CS 80699\\ 62228 Calais Cedex, France\\

\noindent
Claudia Malvenuto\\
Dipartimento di Matematica G. Castelnuovo\\ Sapienza Universit\`a
di Roma\\ P.le Aldo Moro 5\\ 00185, Roma, Italy\\

\noindent
Fr\'ed\'eric Patras\\
LYSM, IRL 2019 CNRS\\permanent address:\\UMR 7351 CNRS\\
        		Universit\'e de Nice\\
        		Parc Valrose\\
        		06108 Nice Cedex 02,
        		France

\end{document}